\documentclass[12pt]{amsart}
\usepackage[T1]{fontenc}
\usepackage[english]{babel}
\usepackage{graphicx}
\usepackage{amsmath,amsfonts}
\usepackage{amssymb}
\usepackage{amsthm}
\usepackage[vlined,ruled,algo2e]{algorithm2e}
\usepackage{tikz}
\usepackage{enumerate}
\usetikzlibrary{arrows,plotmarks,shapes,snakes}
\usepackage{color}

\theoremstyle{plain}
\newtheorem{theorem}{Theorem}
\newtheorem{prop}[theorem]{Proposition}
\newtheorem{cor}[theorem]{Corollary}
\newtheorem{lem}{Lemma}
\newtheorem{exmp}{Example}

\theoremstyle{definition}
\newtheorem{defn}{Definition}

\newcommand{\setpr}[1]{{\mathbb{P}_{#1}}}

\def\P{{\mathcal{P}}}

\definecolor{light-gray}{gray}{0.65}

\begin{document}



\title[Unifying Uncertainty Representations: Clouds]{Unifying practical uncertainty representations: II. Clouds}


\author[S. Destercke]{S\'ebastien Destercke}
\address{Institut de Radioprotection et S\^uret\'e nucl\'eaire, B\^at 720, 13115 St-Paul lez Durance, FRANCE}
\email{sdestercke@gmail.com}
\author[D. Dubois]{Didier Dubois}
\address{Universit\'e Paul Sabatier, IRIT, 118 Route de Narbonne, 31062 Toulouse} \email{dubois@irit.fr}
\author[E. Chojnacki]{Eric Chojnacki}
\address{Institut de Radioprotection et S\^uret\'e nucl\'eaire, B\^at 720, 13115 St-Paul lez Durance, FRANCE} \email{eric.chojnacki@irsn.fr}

\begin{abstract}
There exist many simple tools for jointly capturing variability and incomplete information by means of uncertainty
representations. Among them are random sets, possibility
distributions, probability intervals, and the more recent Ferson's
p-boxes and Neumaier's clouds, both defined by pairs of possibility distributions. In the companion
paper, we have extensively studied a generalized form of p-box and situated it with respect to other models . This paper focuses on the links between
clouds and other representations. Generalized p-boxes are shown to be clouds with comonotonic distributions. In general, clouds
cannot always be represented  by
random sets, in fact not even by $2$-monotone (convex) capacities.
\end{abstract}


\maketitle

\section{Introduction}

There exist many different tools for representing imprecise
probabilities. Usually, the more general, the more difficult they
are to handle. Simpler representations, although less expressive,
usually have the advantage of being more tractable. Over the years,
several such representations have been proposed. Among
them are possibility distributions~\cite{Zadeh78}, probability
intervals~\cite{CamposAll94}, and more recently
p-boxes~\cite{FersonAll03} and
clouds~\cite{Neumaier04,NeumaierWeb04}. Comparing their
respective expressive power is a natural task. Finding formal relations between such
representations also facilitates a unified handling of uncertainty.

\begin{figure}
\begin{center}
\begin{tikzpicture}[scale=0.9]
 \node (ImpProb) at (0,-0.6) [draw] {\scriptsize{Credal sets}};
 \node (LowUpp) at (0,-1.6) [draw] {\scriptsize{Coherent lower/upper probabilities}};
 \node (2mon) at (0,-2.6) [draw] {\scriptsize{2-monotone capacities}} ;
 \node (randset) at (0,-3.8) [draw] {\scriptsize{Random sets ($\infty$-monotone)}} ;
 \node (GenPbox) at (0,-5.6) [draw,ultra thick] {\scriptsize{Generalized p-boxes}} ;
 \node (Pbox) at (0,-6.8) [draw] {\scriptsize{P-boxes}} ;
 \node (Prob) at (0,-8) [draw] {\scriptsize{Probabilities}} ;
 \node (ProbInt) at (-3.3,-4.6) [draw] {\scriptsize{Probability
 Intervals}} ;
 \node (Poss) at (3,-7.8) [draw] {\scriptsize{Possibilities}} ;

 \draw[->] (ImpProb) -- (LowUpp)  ;
 \draw[->]  (LowUpp) -- (2mon) ;
 \draw[->]  (2mon) -- (randset)  ;
 \draw[->,ultra thick]  (randset) -- (GenPbox) ;
 \draw[->,ultra thick] (GenPbox) -- (Pbox)  ;
 \draw[->]  (Pbox) -- (Prob) ;
 \draw[->,ultra thick] (GenPbox.south east) -- (Poss) ;
 \draw[->] (2mon.south west) -- (ProbInt)  ;
 \draw[->]  (ProbInt) -- (Prob.north west) ;
 \draw[dashed,->,ultra thick]  (Poss) .. controls
+(-0.2,1.2) .. (GenPbox.east) ; \draw[dashed,->,ultra thick]
(GenPbox) .. controls +(-0.4,1.2) .. (ProbInt) ;

\end{tikzpicture}  \caption{Relationships among representations. $A \longrightarrow B$ A generalizes B. $A \dashrightarrow B$: B is representable by A }\label{fig:summarypre}
\end{center}
\end{figure}
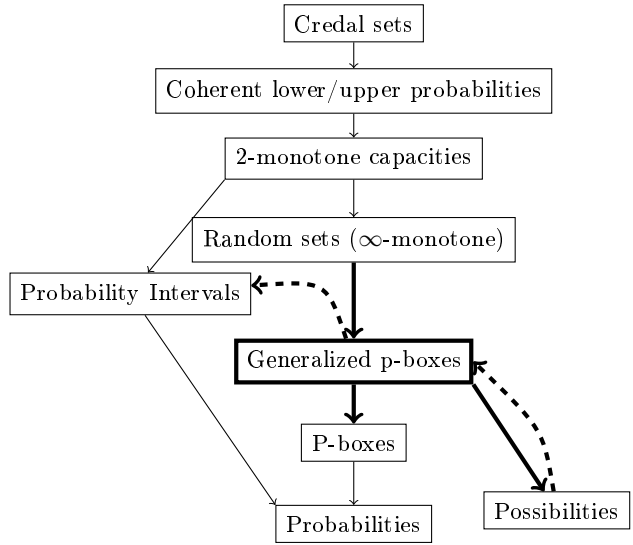

In the first part of paper \cite{DesterckeAll07IJAR}, a generalized notion of p-boxes is studied
and related to representations mentioned above. It is shown that any
generalized p-box is representable by a pair of possibility
distributions, and that generalized p-boxes are special cases of
random sets. Their interpretation in terms of lower and upper confidence bounds on a
collection of nested subsets makes them intuitive  simple
representations. Figure~\ref{fig:summarypre} recalls the connections established in the companion paper between the studied representations, going from the most
(top) to the least (bottom) general.

The present paper completes
Figure~\ref{fig:summarypre} by adding clouds to it, making one step
further towards the unification of uncertainty models. Clouds, encoded by a pair of fuzzy sets, were recently introduced by Neumaier~\cite{Neumaier04} as  a means to cope with imprecision while remaining
computationally tractable, even in high dimensional spaces. Recently,
Fuchs and Neumaier~\cite{FuchsNeum08} have applied clouds to space shuttle design problem,
demonstrating some of the potential of the representation. Moreover,
as clouds are syntactically equivalent to interval-valued fuzzy sets with
some boundary conditions, analyzing their connection with respect to
other uncertainty theories also provides some insight about
how interval-valued fuzzy sets can be interpreted by such theories.
As we shall see, generalized p-boxes, studied in the companion
paper, constitute a bridge between clouds, possibility distributions
and usual p-boxes.

The paper is divided into five main sections;
Section~\ref{sec:Clouds} studies the formalism of clouds and relates
them to pairs of possibility distributions and to generalized
p-boxes. It is shown that generalized p-boxes are equivalent to a
particular subfamily of clouds, named here comonotonic clouds.
Section~\ref{sec:NoncoCloud} studies non-comonotonic clouds. Since the lower
probability they induce are in general even not 2-monotone, simpler outer and inner
approximations are proposed; Section~\ref{sec:probintclouds} then
studies relations between clouds and probability intervals. As
neither of them is a special case of the other, some transformations
of probability intervals into outer-approximating clouds are proposed.
Section~\ref{sec:Cloudrealline} extends some of our results to the
case of continuous models defined on the real line, since such
models are often encountered in applications. The particular case of
thin clouds, for which both upper and lower distributions coincide is emphasized, as they have non-empty credal sets in the
continuous setting.

To make the paper easier to read, longer proofs have been moved to
the appendix. We will often refer to useful results from the  companion paper \cite{DesterckeAll07IJAR}, where basics about other representations
and uncertainty theories considered here can be found. Some definitions are recalled in footnotes. In the first four sections, we
consider that our uncertainty concerns the value that a variable
could assume on a finite set $X$ containing $n$ elements.

\section{Clouds}
\label{sec:Clouds}
Clouds were introduced by Neumaier~\cite{Neumaier04} as a probabilistic generalizations of intervals. 
\begin{defn}
\label{def:clouds} A cloud $[\delta,\pi]$ is defined as a pair of
mappings $\delta: X \to [0,1]$ and $\pi: X \to [0,1]$ from the set
$X$ to the unit interval $[0,1]$, such that \begin{itemize}
\item $\delta$ is pointwise less than or equal to $\pi$ (i.e. $\delta \leq \pi$), 
\item $\pi(x)=1$ for at least one element $x$ in $X$,
\item$\delta(y)=0$ for at least one element $y$ in $X$. 
\end{itemize}
$\delta$ and $\pi$ are respectively the lower and upper
distributions of a cloud.
\end{defn}
As mappings $\delta,\pi$ are
mathematically equivalent to two nested fuzzy membership functions,  a
cloud $[\delta,\pi]$ is mathematically equivalent to an
interval-valued fuzzy set (IVF)\cite{Zadeh75} with boundary conditions
($\pi(x)=1$ and $\delta(y)=0$). More precisely, it is mathematically
equivalent to an interval-valued membership function whereby the
membership value of each element $x$ of $X$ lies in $[\delta(x),\pi(x)]$.
Since a cloud is equivalent to a pair of fuzzy membership functions,
at most $2|X|-2$ values (notwithstanding boundary constraints on
$\delta$ and $\pi$) are needed to fully determine a cloud on a
finite set. Two subcases of clouds considered by
Neumaier~\cite{Neumaier04} are the thin and fuzzy clouds. A
\emph{thin cloud} is defined as a cloud for which $\delta=\pi$,
while a \emph{fuzzy cloud} is a cloud for which $\delta=0$.

 Neumaier defines the credal
set\footnote{A credal set $\P$ is a closed convex set of probability distributions, here described by constraints on
probabilities of some events. 
} 
$\P_{[\delta,\pi]}$ induced by a cloud $[\delta,\pi]$, as:
\begin{equation} \P_{[\delta,\pi]}\!=\!\{P \in \setpr{X} | P( \{x \in X|\delta(x)
\geq \alpha \} ) \leq 1- \alpha \leq P( \{x \in X|\pi(x) > \alpha \}
) \}
\end{equation}
where $\setpr{X}$ is the set of probability measures on $X$.
Interestingly enough, this definition gives a mean to interpret IVF
sets in terms of credal sets, or in terms of imprecise
probabilities, eventually ending up with a behavioral interpretation
of IVF by using Walley's~\cite{Walley91} theory of imprecise
probabilities.

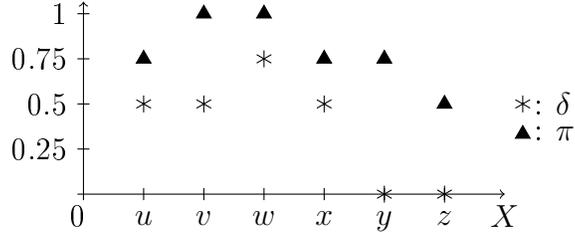
\begin{figure}
\begin{center}
\begin{tikzpicture}[scale=0.8]
\draw[->] (-0.1,0) node[below] {0} -- (7,0) node[below] {$X$} ;
\draw[->] (0,-0.1) -- (0,3.2) ; \draw (0.1,3) -- (-0.1,3) node[left]
{1} ; \draw (1,0.1) --  (1,-0.1) node[below] {$u$}; \draw (2,0.1) --
(2,-0.1) node[below] {$v$}; \draw (3,0.1) -- (3,-0.1) node[below]
{$w$}; \draw (4,0.1) --  (4,-0.1) node[below] {$x$}; \draw (5,0.1)
--  (5,-0.1) node[below] {$y$}; \draw (6,0.1) -- (6,-0.1)
node[below] {$z$}; \draw (0.1,0.75) -- (-0.1,0.75) node[left] {0.25}
; \draw (0.1,1.5) -- (-0.1,1.5) node[left] {0.5} ; \draw (0.1,2.25)
-- (-0.1,2.25) node[left] {0.75} ;

\draw plot[only marks,mark=triangle*,mark options={scale=2.0}]
coordinates {(1,2.25) (2,3) (3,3) (4,2.25) (5,2.25) (6,1.5)} ;

\draw plot[only marks,mark=asterisk,mark options={scale=2.0}]
coordinates {(1,1.5) (2,1.5) (3,2.25) (4,1.5) (5,0) (6,0)};

\draw plot[mark=asterisk,mark options={scale=2.0}] coordinates
{(7.3,1.5)} ; \node[right] at (7.3,1.5) {: $\delta$} ;

\draw plot[mark=triangle*,mark options={scale=2.0}] coordinates
{(7.3,1)} ; \node[right] at (7.3,1) {: $\pi$} ;

\end{tikzpicture}
\caption{Cloud $[\delta,\pi]$ of Example~\ref{exmp:clouds1}}
\label{fig:CloudEx}
\end{center}
\end{figure}

Let $0=\gamma_0 < \gamma_1 < \ldots < \gamma_M =1$ be the ordered
distinct values taken by both $\delta$ and $\pi$ on elements of $X$,
then denote the strong and regular cuts as
\begin{equation}
\label{eq:clouduppercut} B_{\overline{\gamma}_i}= \{ x \in X |
\pi(x)
> \gamma_{i}\} \mbox{~and~}
 B_{\gamma_i} = \{ x \in X | \pi(x) \geq \gamma_{i} \}
 \end{equation}
for the upper distribution $\pi$ and
\begin{equation}
\label{eq:cloudlowercut}C_{\overline{\gamma}_i}= \{ x \in X |
\delta(x)
> \gamma_{i}\} \mbox{~and~}
 C_{\gamma_i} = \{ x \in X | \delta(x) \geq \gamma_{i} \}\end{equation}
for the lower distribution $\delta$. Note that in the finite case,
$B_{\overline{\gamma}_i} = B_{\gamma_{i+1}}$ and
$C_{\overline{\gamma}_i} = C_{\gamma_{i+1}}$, with $\gamma_{M+1}=1$,
and also
$$\emptyset = B_{\overline{\gamma}_{M}} \subset B_{\overline{\gamma}_{M-1}}
\subseteq \ldots
 \subseteq B_{\overline{\gamma}_{0}} = X ; \quad \emptyset =
 C_{\gamma_M} \subseteq C_{\gamma_{M-1}} \subseteq  \ldots \subseteq
 C_{\gamma_{0}}= X$$
 and  since $\delta
\leq \pi$, this implies that $C_{\gamma_{i}} \subseteq
B_{\gamma_{i}}$, hence $C_{\gamma_{i}} \subseteq
B_{\overline{\gamma}_{i-1}}$, for all $i = 1, \dots, M$. In such a
finite case, a cloud is said to be discrete. In terms of constraints
bearing on probabilities, the credal set $\P_{[\delta,\pi]}$ of a
finite cloud is equivalently defined by the finite set of inequalities:
\begin{equation}
\label{eq:cloudcons} i=0,\ldots,M, \quad P(C_{\gamma_i}) \leq 1 -
\gamma_{i} \leq P(B_{\overline{\gamma}_{i}})
\end{equation}
under the above inclusion constraints. Note that some conditions, in
addition to boundary ones advocated in Definition~\ref{def:clouds},
must hold for $\P_{[\delta,\pi]}$ to be non-empty in the finite
case. In particular, distribution $\delta$ cannot be equal to
$\pi$ (i.e. $\delta(x) \neq \pi$). Otherwise, we have $C_{\gamma_{i}} =
B_{\overline{\gamma}_{i-1}} (= B_{\gamma_{i}} )$, that is $\pi$ and
$\delta$ have a common $\gamma_i$-cut, and there
 is no probability distribution satisfying the constraint $1 -
\gamma_{i-1} \leq P(C_{\gamma_{i}}) \leq 1 - \gamma_{i}$ since $
\gamma_{i-1} < \gamma_i$. So, thin finite clouds induce empty credal
sets.
\begin{exmp}
\label{exmp:clouds1} This example illustrates the notion of a cloud
and will be used in the next sections to illustrate various results.
Let us consider a set $X=\{u,v,w,x,y,z\}$ and the following cloud
$[\delta,\pi]$, pictured in Figure~\ref{fig:CloudEx}, defined on
this set:
\begin{displaymath}
\begin{array}{c@{\hspace{10mm}}c@{\hspace{5mm}}c@{\hspace{5mm}}c@{\hspace{5mm}}c@{\hspace{5mm}}c@{\hspace{5mm}}c}
& u & v & w & x & y & z \\
\hline
\pi & 0.75 & 1 & 1 & 0.75 & 0.75 & 0.5 \\
\delta & 0.5 & 0.5 & 0.75 & 0.5 & 0 & 0 \\
\hline
\end{array}
\end{displaymath}
 The values
$\gamma_i$ corresponding to this cloud are
\begin{align*}
0 \leq 0.5 \leq 0.75 \leq 1 \\
\gamma_0 \leq \gamma_1 \leq \gamma_2 \leq \gamma_3.
\end{align*}
Constraints associated to this cloud and corresponding to
Equation (\ref{eq:cloudcons}) are
\begin{align*}
P(C_{\gamma_3} = \emptyset) & \leq 1 - 1 \leq
P(B_{\overline{\gamma}_3} = \emptyset) \\
P(C_{\gamma_2} = \{w\}) & \leq 1 - 0.75 \leq  P(B_{\overline{\gamma}_2} = \{v,w\}) \\
P(C_{\gamma_1} = \{u,v,w,x \}) & \leq 1 - 0.5 \leq
P(B_{\overline{\gamma}_1} = \{u,v,w,x,y\}) \\
P(C_{\gamma_0} = X) & \leq 1 - 0 \leq P(B_{\overline{\gamma}_0} = X)
\end{align*}
\end{exmp}

\subsection{Clouds in the setting of  possibility theory}
\label{sec:CloudPoss}

To relate clouds with possibility distributions\footnote{A
possibility distribution is a mapping $\pi: X \to [0,1]$, with
$\pi(x)=1$ for at least one element, and inducing a credal set
$\P_\pi$ such that $P \in \P_\pi$ iff $1- \alpha \leq P( \{ x \in X
| \pi(x)
> \alpha \}$ for all $\alpha \in [0,1]$}, first consider the case of
\emph{fuzzy clouds} $[\delta,\pi]$. In this case, $\delta=0$ and
$C_{\gamma_i}=\emptyset$ for $i=1,\ldots,M$, which means that
constraints given by Equations (\ref{eq:cloudcons}) reduce to $ 1 -
\gamma_{i} \leq P(B_{\overline{\gamma}_{i}})$ for $i=0,\ldots,M$
which, by using Proposition 2.5 of the  companion paper \cite{DesterckeAll07IJAR}, induces a
credal set $\P_{\pi}$ equivalent to the one induced by the
possibility distribution $\pi$. This shows that \emph{fuzzy clouds}
are equivalent to possibility distributions. The following
proposition is a direct consequence of this observation:
\begin{prop}
\label{prop:CloudPoss} Uncertainty modeled by a cloud $[\delta,\pi]$
is representable by the pair of possibility distributions $1-\delta$
and $\pi$, and we have: $$\P_{[\delta,\pi]}= \P_{\pi} \cap
\P_{1-\delta}$$
\end{prop}
\begin{proof}[\textbf{Proof of Proposition~\ref{prop:CloudPoss}}]
Consider a cloud $[\delta, \pi]$ and the constraints inducing the
credal set $\P_{[\delta,\pi]}$. As for generalized p-boxes, these
constraints can be split into two sets of constraints, namely, for
$i=0,\ldots,M$, $P(C_{\gamma_i}) \leq 1 - \gamma_{i}$ and $1 -
\gamma_{i} \leq P(B_{\overline{\gamma}_{i}})$. Since
$B_{\overline{\gamma}_{i}}$ are strong cuts of $\pi$, then by
Proposition 2.5. in \cite{DesterckeAll07IJAR} we know that these constraints
define a credal set equivalent to $\P_{\pi}$.

Note then that $P(C_{\gamma_i}) \leq 1 - \gamma_{i}$ is equivalent
to \mbox{$P(C_{\gamma_i}^c) \geq \gamma_i$} (where
$C_{\gamma_i}^c=\{ x \in X | 1 - \delta(x) > 1 - \gamma_{i} \}$). By
construction, $1-\delta$ is a normalized possibility distribution.
Interpreting these inequalities in the light of Proposition 2.5.
in \cite{DesterckeAll07IJAR}, we see that they define the credal set
$\P_{1-\delta}$. By merging the two set of constraints, we get
$\P_{\delta,\pi}= \P_{\pi} \cap \P_{1-\delta}$.
\end{proof}

This proposition shows that, as for generalized p-boxes, the credal
set induced by a cloud is representable by a pair of possibility
distributions~\cite{DuboisPrade05}. This analogy between generalized
p-boxes and clouds is studied in
Section~\ref{sec:CloudGenPbox}. This result also confirms that a
cloud $[\delta, \pi]$ is equivalent to its {\em mirror} cloud $[1-\pi, 1-
\delta]$ ($1-\pi$ becoming the lower distribution, and $1-
\delta$ the upper one).
\begin{exmp}
\label{exmp:clouds2} Possibility distributions $\pi,1-\delta$
representing the cloud of Example~\ref{exmp:clouds1} are:
\begin{displaymath}
\begin{array}{c@{\hspace{10mm}}c@{\hspace{5mm}}c@{\hspace{5mm}}c@{\hspace{5mm}}c@{\hspace{5mm}}c@{\hspace{5mm}}c}
& u & v & w & x & y & z \\
\hline
\pi & 0.75 & 1 & 1 & 0.75 & 0.75 & 0.5 \\
1 - \delta & 0.5 & 0.5 & 0.25 & 0.5 & 1 & 1 \\
\hline
\end{array}
\end{displaymath}
\end{exmp}

\subsection{Clouds with non-empty credal sets} \label{Nonemptycredal}

We now explore under which conditions a cloud $[\delta,\pi]$ induces
a non-empty credal set $\P_{[\delta,\pi]}$. Using
the fact that clouds are representable by pairs of possibility
distributions, and applying Chateauneuf~\cite{Chateauneuf94}
characteristic condition ($\forall A \subset X, Bel_1(A) + Bel_2(A^c) \leq 1$) under which the credal sets associated to
two belief functions $ Bel_1$ and $Bel_2$ have a non-empty intersection, the following
necessary and sufficient condition obtains:
\begin{prop} A cloud  $[ \delta, \pi]$ has a non-empty credal set if and only if
$$\forall A \subseteq X, \max_{x\in A} \pi(x) \geq \min_{y\not\in A} \delta(y) $$
\label{Chateau}
\end{prop}
\begin{proof} Chateauneuf's condition applied to possibility distributions $\pi_1$ and $\pi_2$ reads $\forall A \subseteq X, \Pi_1(A)+\Pi_2(A^c)\geq 1$. Choose $\pi_1 = \pi$ and $\pi_2 = 1-\delta$. In particular $\Pi_2(A^c) = 1 -  \min_{y\not\in A} \delta(y)$.
\end{proof}

A naive test for non-emptiness based on Proposition~\ref{Chateau}
would have exponential complexity, but in the case of clouds, it can
be simplified as follows: suppose the space $X=\{x_1,\ldots,x_n\}$
is indexed such that $\pi(x_1) \leq \pi(x_2)\dots \leq \pi(x_n)=1$
and consider an event $A$ such that $\max_{x\in A} \pi(x) =
\pi(x_i)$. The tightest constraint of the form $\max_{x\in A} \pi(x)
=\pi(x_i) \geq \min_{y\not\in A} \delta(y)$ is when choosing $A =
\{x_1, \dots x_i\}$. Checking non-emptiness then comes down to
checking the following set of $n-1$ inequalities:
\begin{equation} \label{eq:consnonempt} j = 1,\ldots, n-1 \quad \pi(x_i)
\geq \min_{j > i} \delta(x_j).
\end{equation}
This gives us an efficient tool to check the non-emptiness of a
given cloud on a finite set, or to build a non-empty cloud from the knowledge of
either $\delta$ or $\pi$. For instance, knowing $\delta$, the cloud
$[\delta, \pi]$ such that $\pi(x_i) = \min_{j > i} \delta(x_j),  j =
1,\ldots, n-1$ is the most restrictive non-empty cloud one may
build, assuming the ordering  $\pi(x_1) \leq \pi(x_2)\dots \leq
\pi(x_n)=1$ (changing this assumption yields another non-empty
cloud).

Now, consider the extreme case of a cloud for which $C_{\gamma_i} =
B_{\overline{\gamma}_{i}}$ for all $i=1,\ldots,M$ in equation
(\ref{eq:cloudcons}). In this case, $P(B_{\overline{\gamma}_{i}}) =
P(C_{\gamma_i}) = 1 - \gamma_{i}$ for all $i=1,\ldots,M$. Suppose
distribution $\pi$ takes distinct values on all elements of $X$.
Ordering elements of $X$ by increasing values of $\pi(x)$ ($\forall
i, \pi(x_i)
> \pi(x_{i-1})$) enforces $\delta(x_i) = \pi(x_{i-1})$, with
$\delta(x_1) = 0$. Let $\delta_\pi$ be this lower distribution. The
(almost thin) cloud $[\delta_\pi, \pi]$ satisfies equations
(\ref{eq:consnonempt}), and since $P(B_{\overline{\gamma}_{i}}) = 1
- \gamma_{i}$, the induced credal set $\P_{[\delta_\pi,\pi]}$
contains the single probability measure $P$ with distribution
$p(x_i)=\pi(x_i)-\pi(x_{i-1})$ for all $x_i \in X$, with
$\pi(x_0)=0$. So if a finite cloud $[ \delta, \pi]$  is such that if
$\delta > \delta_\pi$, it induces an empty credal set
$\P_{[\delta,\pi]}$; and if $\delta \leq \delta_\pi$, then the
induced credal set $\P_{[\delta,\pi]}$ is not empty.

Equations (\ref{eq:consnonempt}) can be extended to the case
of any two possibility distributions $\pi_1,\pi_2$ for which we want
to check whether $\P_{\pi_1} \cap \P_{\pi_2}$ is empty or not. This is meaningful because the setting of clouds does not
cover all pairs  $\pi_1,\pi_2$ such that  $\P_{\pi_1} \cap \P_{\pi_2}\neq\emptyset$.
To check it, first recall that for any two possibility
distributions $\pi_1,\pi_2$, we do have $\P_{\min(\pi_1,\pi_2)}
\subseteq \P_{\pi_1} \cap \P_{\pi_2}$, but, in general, the
converse inclusion \cite{DuboisPrade90b} does not hold. From this remark, we have
\begin{itemize}
\item  $ \P_{\pi_1} \cap \P_{\pi_2} \neq \emptyset$ as soon as $\min(\pi_1,\pi_2)$ is a normalized possibility distribution.
\item Not all pairs $\pi_1,\pi_2$ such that $ \P_{\pi_1} \cap \P_{\pi_2} \neq \emptyset$ derive
 from a cloud $[ 1 - \pi_2, \pi_1 ]$. Indeed, the normalization of $\min(\pi_1,\pi_2)$ does not imply  that $1 - \pi_2 \leq \pi_1$.
\end{itemize}

\subsection{Generalized p-boxes as a special kind of clouds}
\label{sec:CloudGenPbox}

We remind~ \cite{DesterckeAll07IJAR}  that a generalized p-box $[\underline{F},\overline{F}]$ is
defined by two comonotonic mappings $\underline{F}:X\to[0,1]$,
$\overline{F}:X\to[0,1]$ with $\underline{F}\leq\overline{F}$ and
$\underline{F}(x)=\overline{F}(x)=1$ for at least one element $x$ of
$X$. They induce a pre-order $\leq_{[\underline{F},\overline{F}]}$
on $X$ such that $x \leq_{[\underline{F},\overline{F}]} y$ if
$\underline{F}(x) \leq \underline{F}(y)$ and  $\overline{F}(x) \leq
\overline{F}(y)$, and elements of $X$ are here indexed such that $i
\leq j$ implies $x_i \leq_{[\underline{F},\overline{F}]} x_j$. A
generalized p-box $[\underline{F},\overline{F}]$ induces the
following credal set:
\begin{equation}\label{eq:PBoxCons}\P_{[\underline{F},\overline{F}]}=\{P \in \setpr{X}
|i=1,\ldots,n, \; \alpha_i \leq P(A_i) \leq \beta_i
\}.\end{equation} Where $A_i=\{x_1,\ldots,x_i\}$,
$\alpha_i=\underline{F}(x_i)$ and $\beta_i=\overline{F}(x_i)$ are
lower and upper confidence bounds on set $A_i$. Note that
$A_1\subseteq\ldots\subseteq A_n$, $\alpha_1 \leq \ldots \leq
\alpha_n$ and $\beta_1 \leq \ldots \leq \beta_n$. The proposition
below lays bare the nature of the relationship between such
generalized p-boxes and clouds:
\begin{prop}
\label{lem:CloudBel1} Let $[\delta,\pi]$ be a cloud defined on $X$.
Then, the three following statements are equivalent:
\begin{enumerate}[(i)]
\item \label{cloudbelit1} The cloud $[\delta,\pi]$ can be encoded as a generalized
p-box $[\underline{F},\overline{F}]$ such that
$\P_{[\delta,\pi]}=\P_{[\underline{F},\overline{F}]}$
\item \label{cloudbelit2} $\delta$ and $\pi$ are comonotonic ($\delta(x)<\delta(y) \Rightarrow \pi(x) \leq \pi(y)$)
\item \label{cloudbelit3}  Sets $ \{B_{\overline{\gamma}_i},C_{\gamma_j} |
i,j=0,\ldots ,M \}$ form a nested sequence (i.e. are completely
(pre-)ordered with respect to inclusion).
\end{enumerate}
\end{prop}

\begin{proof}[\textbf{Proof of Proposition~\ref{lem:CloudBel1}}]
We use a cyclic proof to show that statements (\ref{cloudbelit1}),
(\ref{cloudbelit2}), (\ref{cloudbelit3}) are equivalent.

(\ref{cloudbelit1})$\Rightarrow$(\ref{cloudbelit2}) From the assumption, $\delta=1-\pi_{\underline{F}}$ and $\pi=\pi_{\overline{F}}$. Hence, using Proposition 3.3 in \cite{DesterckeAll07IJAR} and
the definition of a generalized p-box, $\delta$ and $\pi$ are comonotone, hence (\ref{cloudbelit1})$\Rightarrow$(\ref{cloudbelit2}). 


(\ref{cloudbelit2})$\Rightarrow$(\ref{cloudbelit3}) we will show
that
 if (\ref{cloudbelit3}) does not hold, then (\ref{cloudbelit2}) does not hold either. Assume sets $ \{B_{\overline{\gamma}_i},C_{\gamma_j}
| i,j=0,\ldots ,M \}$ do not form a nested sequence, meaning that
there exists two sets $C_{\gamma_j}, B_{\overline{\gamma}_i}$ with
$j<i$ s.t. $C_{\gamma_j} \not \subset B_{\overline{\gamma}_i}$ and
$B_{\overline{\gamma}_i} \not \subset C_{\gamma_j}$. This is
equivalent to asserting $\exists x,y \in X$ such that $\delta(x)\geq
\gamma_j$, $\pi(x)\leq \gamma_i$, $\delta(y)<\gamma_j$ and
$\pi(y)>\gamma_i$. This implies $\delta(y)<\delta(x)$ and $\pi(x) <
\pi(y)$, and that $\delta,\pi$ are not comonotonic.

(\ref{cloudbelit3})$\Rightarrow$(\ref{cloudbelit1}) Assume the sets
$B_{\overline{\gamma}_i}$ and $C_{\gamma_j}$ form a globally nested
sequence whose current element is $A_k$.  Then the set of
constraints defining a cloud can be rewritten in the form $\alpha_k
\leq P(A_k) \leq \beta_k$, where  $\alpha_k = 1 - \gamma_{i}$ and
$\beta_k = \min\{1- \gamma_{j}| B_{\overline{\gamma}_i} \subseteq
C_{\gamma_j} \}$ if $A_k = B_{\overline{\gamma}_i}$; $\beta_k = 1 -
\gamma_{i}$ and \mbox{$\alpha_k = \max\{1- \gamma_{j}|
B_{\overline{\gamma_j}} \subseteq C_{\gamma_i}\}$} if $A_k =
C_{\gamma_i}$. Since $0 = \gamma_0 < \alpha_1 < \ldots < \alpha_M
=1$, these constraints are equivalent to those describing a
generalized p-box (Equations (\ref{eq:PBoxCons})). This ends the
proof.
\end{proof}

Proposition~\ref{lem:CloudBel1} indicates that only those clouds for
which $\delta$ and $\pi$ are comonotonic can be encoded by
generalized p-boxes, and from now on, we shall call such clouds
\emph{comonotonic}. Using proposition 3.3 of the  companion paper \cite{DesterckeAll07IJAR} and
given a comonotonic cloud $[\delta,\pi]$, we can express this cloud
as the following generalized p-box $\underline{F},\overline{F}$
defined for any $x \in X$:
\begin{equation}
\label{eq:cloudtogenpbox} \overline{F}(x)=\pi(x) \textrm{ and }
\underline{F}(x)=\min\{\delta(y) | y \in X, \delta(y) > \delta(x)\}.
\end{equation}
Conversely, note that any generalized p-box
$[\underline{F},\overline{F}]$ can be encoded by a comonotonic
cloud, simply taking $\delta=1-\pi_{\underline{F}}$ and
$\pi=\pi_{\overline{F}}$ (see Proposition 3.3 in \cite{DesterckeAll07IJAR}).
This means that generalized p-boxes are special cases of clouds, and that comonotonic clouds and generalized
p-boxes are equivalent representations. Also note that a comonotonic
cloud $[\delta,\pi]$ and the equivalent generalized p-box
$[\underline{F},\overline{F}]$ induce the same complete pre-order
on elements of $X$, that we note
$\leq_{[\underline{F},\overline{F}]}$ to remain coherent with the
notations of the  companion paper \cite{DesterckeAll07IJAR}. We consider that elements $x$ of
$X$ are indexed accordingly, as already specified.

In practice, this means that all the results that hold for
generalized p-boxes also hold for comonotonic clouds, and
conversely. In particular, comonotonic clouds are special cases of
random sets\footnote{A random set is a non-negative mapping
$m:\wp(X)\to[0,1]$ such that $\sum_{E \subset X} m(E)=1; m(\emptyset) = 0$. It is also
completely characterized by the belief function $Bel$ such that $ \forall A \subset X$, $Bel(A)=\sum_{E \subseteq
A}m(E)$. The credal set $\P_{Bel}$ induced by such a random set is $\P_{Bel}= \{ P \in \setpr{X} | \forall A \subseteq X,
Bel(A) \leq P(A) \}$}, in the sense that, for any comonotonic cloud
$[\delta,\pi]$, there is a belief function $Bel$ such that
$\P_{[\delta,\pi]}=\P_{Bel}$. Adapting Equations (13) of the
 companion paper \cite{DesterckeAll07IJAR} to the case of a comonotonic cloud $[\delta,\pi]$,
this random set is such that, for $j=1,\ldots,M$:
\begin{equation}
\label{eq:cloudRStransform} \left\{ \begin{array}{c} E_j=\{ x \in X
| (\pi(x) \geq \gamma_j)
\wedge (\delta(x) < \gamma_j) \} \\
m(E_j)=\gamma_j - \gamma_{j-1}.
\end{array} \right.
\end{equation}
Note that in the formalism of clouds this random set can be
expressed in terms of the sets $
\{B_{\overline{\gamma_i}},C_{\gamma_i} | i=0,\ldots ,M \}$. Namely,
for $j=1,\ldots,M$:
\begin{equation}
\label{eq:cloudcutRStransform} \left\{ \begin{array}{c} E_j=
B_{\overline{\gamma}_{j-1}} \setminus C_{\gamma_j} = B_{\gamma_j} \setminus C_{\gamma_j} \\
m(E_j)=\gamma_j - \gamma_{j-1}.
\end{array} \right.
\end{equation}

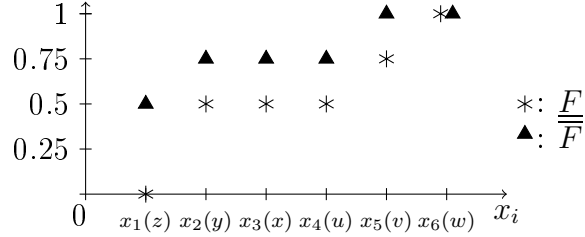
\begin{figure}
\begin{center}
\begin{tikzpicture}[scale=0.8]
\draw[->] (-0.1,0) node[below] {0} -- (7,0) node[below] {$x_i$} ;
\draw[->] (0,-0.1) -- (0,3.2) ; \draw (0.1,3) -- (-0.1,3) node[left]
{1} ; \draw (1,0.1) --  (1,-0.1) node[below] {$\scriptstyle
x_1(z)$}; \draw (2,0.1) -- (2,-0.1) node[below] {$\scriptstyle
x_2(y)$}; \draw (3,0.1) -- (3,-0.1) node[below] {$\scriptstyle
x_3(x)$}; \draw (4,0.1) -- (4,-0.1) node[below] {$\scriptstyle
x_4(u)$}; \draw (5,0.1) -- (5,-0.1) node[below] {$\scriptstyle
x_5(v)$}; \draw (6,0.1) -- (6,-0.1) node[below] {$\scriptstyle
x_6(w)$}; \draw (0.1,0.75) -- (-0.1,0.75) node[left] {0.25} ; \draw
(0.1,1.5) -- (-0.1,1.5) node[left] {0.5} ; \draw (0.1,2.25) --
(-0.1,2.25) node[left] {0.75} ; \draw (0.1,3) -- (-0.1,3) node[left]
{1} ;

\draw plot[only marks,mark=triangle*,mark options={scale=2.0}]
coordinates {(1,1.5) (2,2.25) (3,2.25) (4,2.25) (5,3) (6.1,3)} ;

\draw plot[only marks,mark=asterisk,mark options={scale=2.0}]
coordinates {(1,0) (2,1.5) (3,1.5) (4,1.5) (5,2.25) (5.9,3)};

\draw plot[mark=asterisk,mark options={scale=2.0}] coordinates
{(7.3,1.5)} ; \node[right] at (7.3,1.5) {: $\underline{F}$} ;

\draw plot[mark=triangle*,mark options={scale=2.0}] coordinates
{(7.3,1)} ; \node[right] at (7.3,1) {: $\overline{F}$} ;

\end{tikzpicture}
\caption{Generalized p-box $[\underline{F},\overline{F}]$
corresponding to cloud of Example~\ref{exmp:clouds1}}.
\label{fig:CloudTopboxEx}
\end{center}
\end{figure}

\begin{exmp}
\label{exmp:clouds3} From the cloud in Example~\ref{exmp:clouds1},
 $C_{\gamma_3} \subset C_{\gamma_2} \subset
B_{\overline{\gamma_2}} \subset C_{\gamma_1} \subset
B_{\overline{\gamma_1}} \subset B_{\overline{\gamma_0}}$, and the
constraints defining $\P_{[\delta,\pi]}$ can be transformed into
\begin{align*}
0 \leq & C_{\gamma_2}= \{w\}\leq 0.25 \\
0.25 \leq & B_{\overline{\gamma_2}}= \{v,w\}  \leq 0.5 \\
0.25 \leq & C_{\gamma_1} = \{u,v,w,x \} \leq 0.5 \\
0.5 \leq & B_{\overline{\gamma_1}} = \{u,v,w,x,y\}  \leq 1.
\end{align*}
They are equivalent to the generalized p-box
$[\underline{F},\overline{F}]$ pictured on
Figure~\ref{fig:CloudTopboxEx}:
\begin{displaymath}
\begin{array}{c@{\hspace{10mm}}c@{\hspace{5mm}}c@{\hspace{5mm}}c@{\hspace{5mm}}c@{\hspace{5mm}}c@{\hspace{5mm}}c}
& u & v & w & x & y & z \\
\hline
\overline{F} & 0.75 & 1 & 1 & 0.75 & 0.75 & 0.5 \\
\underline{F} & 0.5 & 0.75 & 1 & 0.5 & 0.5 & 0 \\ \hline
\end{array}
\end{displaymath}
The following ordering is compatible with the two distributions (see
Figure~\ref{fig:CloudTopboxEx}): $
 z <_{[\underline{F},\overline{F}]}  y
 <_{[\underline{F},\overline{F}]} x =_{[\underline{F},\overline{F}]} u
 <_{[\underline{F},\overline{F}]} v  <_{[\underline{F},\overline{F}]} w
$

And the corresponding random set, given by Equations
(\ref{eq:cloudcutRStransform}) or (\ref{eq:cloudRStransform}), is:
\begin{displaymath}
m(\{x_5,x_6\}) = 0.25  ;\;  m(\{x_2,x_3,x_4,x_5\}) = 0.25  ;\;
m(\{x_1,x_2\}) = 0.5
\end{displaymath}
\end{exmp}

Comonotonic clouds being special cases of clouds, we may wonder if
some of the results presented in this section extend to clouds that
are not comonotonic (and called non-comonotonic). In particular, can
uncertainty modeled by a non-comonotonic cloud be exactly modeled by
an equivalent random set?


\section{The Nature  of Non-comonotonic Clouds}
\label{sec:NoncoCloud}

When $[\delta,\pi]$ is a non-comonotonic cloud,
Proposition~\ref{prop:CloudPoss} linking clouds and possibility
distributions still holds, but Proposition~\ref{lem:CloudBel1} does
not hold any longer. As we shall see, non-comonotonic clouds appear
to be less interesting, at least from a practical point of view,
than comonotonic ones.

\subsection{Characterization}
\label{sec:NoncoCloudCharac}

\begin{figure}
\begin{center}
\begin{tikzpicture}[scale=0.8]
\draw[->] (-0.1,0) node[below] {0} -- (7,0) node[below] {$X$} ;
\draw[->] (0,-0.1) -- (0,3.2) ; \draw (2,0.1) -- (2,-0.1)
node[below] {$v$}; \draw (3,0.1) -- (3,-0.1) node[below] {$w$};
\draw (4,0.1) --  (4,-0.1) node[below] {$x$}; \draw (5,0.1) --
(5,-0.1) node[below] {$y$}; \draw (6,0.1) -- (6,-0.1) node[below]
{$z$}; \draw (0.1,0.75) -- (-0.1,0.75) node[left] {0.25} ; \draw
(0.1,1.5) -- (-0.1,1.5) node[left] {0.5} ; \draw (0.1,2.25) --
(-0.1,2.25) node[left] {0.75} ; draw (0.1,3) -- (-0.1,3) node[left]
{1} ;

\draw plot[only marks,mark=triangle*,mark options={scale=2.0}]
coordinates { (2,3) (3,3) (4,1.5) (5,1.5) (6,0.75)} ;

\draw plot[only marks,mark=asterisk,mark options={scale=2.0}]
coordinates { (2,0) (3,1.5) (4,0.75) (5,0) (6,0)};

\draw plot[mark=asterisk,mark options={scale=2.0}] coordinates
{(7.3,1.5)} ; \node[right] at (7.3,1.5) {: $\delta$} ;

\draw plot[mark=triangle*,mark options={scale=2.0}] coordinates
{(7.3,1)} ; \node[right] at (7.3,1) {: $\pi$} ;

\end{tikzpicture}
\caption{Cloud $[\delta,\pi]$ of Example~\ref{exmp:noncoclouds}}
\label{fig:NoncoClouds}
\end{center}
\end{figure}
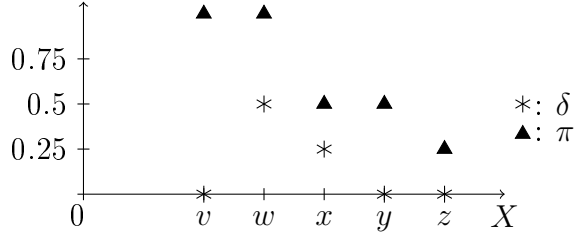

One way of characterizing an uncertainty model is to find the
maximal natural number $n$ such that the lower probability\footnote{The lower probability $\underline{P}$
induced by a credal set $\P$ is $\underline{P}(A)=\min_{P \in \P}P(A)$ for any $A
\subseteq X$.} induced by
this uncertainty model is always $n$-monotone (see \cite{DesterckeAll07IJAR}
or Chateauneuf and Jaffray~\cite{ChateauneufJaffray91} for further details on $n$-monotonicity\footnote{Here we only need 2-monotonicity: A set-function $g$ with domain $2^X$ is 2-monotone if and only if $\forall A, B \subseteq X, g(A) + g(B) \leq g(A\cup B) + g(A\cap B)$.}). This is how we will proceed with non-comonotonic
clouds: let $[\delta,\pi]$ be a non-comonotonic cloud, and
$\P_{[\delta,\pi]}$ the induced credal set. The question is:  what
is the (minimal) $n$-monotonicity of the associated lower
probability $\underline{P}$ induced by $\P_{[\delta,\pi]}$? To address this
question, let us start with an example:

\begin{exmp}
\label{exmp:noncoclouds} Consider a set $X$ with five elements
$\{v,w,x,y,z\}$ and the following non-comonotonic cloud
$[\delta,\pi]$ pictured on Figure~\ref{fig:NoncoClouds}:
\begin{displaymath}
\begin{array}{c@{\hspace{10mm}}c@{\hspace{5mm}}c@{\hspace{5mm}}c@{\hspace{5mm}}c@{\hspace{5mm}}c}
 & v & w & x & y & z \\
 \hline
\pi & 1 & 1 & 0.5 & 0.5 & 0.25  \\
\delta & 0 & 0.5 & 0.25 & 0 & 0  \\
\hline
\end{array}
\end{displaymath}
This cloud is non-comonotonic, since $\pi(v)>\pi(x)$ and
$\delta(v)<\delta(x)$. The credal set $\P_{[\delta,\pi]}$ can also
be defined by the following constraints:
\begin{align*}
P(C_{\gamma_2} = \{w\}) & \leq 1 - 0.5 \leq  P(B_{\overline{\gamma}_2} = \{v,w\}) \\
P(C_{\gamma_1} = \{w,x \}) & \leq 1 - 0.25 \leq
P(B_{\overline{\gamma}_1} = \{v,w,x,y\})
\end{align*}
with $\gamma_2 = 0.5$ and $\gamma_1 = 0.25$. Now, consider the
events
$B_{\overline{\gamma}_2},C_{\gamma_1}^c,B_{\overline{\gamma}_2} \cap
C_{\gamma_1}^c$, $B_{\overline{\gamma}_2} \cup C_{\gamma_1}^c$. We
can check that
\begin{align*}
\underline{P}(B_{\overline{\gamma}_2})= 0.5 & & \underline{P}(C_{\gamma_1}^c)= 0.25 \\
\underline{P}(B_{\overline{\gamma}_2} \cap C_{\gamma_1}^c = \{ v\})=
0 & & \underline{P}(B_{\overline{\gamma}_2} \cup C_{\gamma_1}^c = \{
v,w,y,z \}) = 0.5
\end{align*}
since at most a 0.5 probability mass can be assigned to $x$. Then
the inequality $\underline{P}(B_{\overline{\gamma}_2} \cap
C_{\gamma_1}^c) + \underline{P}(B_{\overline{\gamma}_2} \cup
C_{\gamma_1}^c) < \underline{P}(B_{\overline{\gamma}_2}) +
\underline{P}(C_{\gamma_1}^c)$ holds, indicating that the lower
probability induced by the cloud is not 2-monotone.
\end{exmp}

This example shows that at least some
non-comonotonic clouds induce lower probability measures that are
not $2$-monotone. The following proposition gives a
general characterization of a large family of such non-comonotonic clouds:
\begin{prop}
\label{prop:CloudBel} Let $[\delta, \pi]$ be a non-comonotonic cloud
and assume there is a pair of events
$B_{\overline{\gamma}_i},C_{\gamma_j}$ in the cloud s.t.
$B_{\overline{\gamma}_i} \cap C_{\gamma_j} \not\in
\{B_{\overline{\gamma}_i},C_{\gamma_j},\emptyset \}$ and
$B_{\overline{\gamma}_i} \cup C_{\gamma_j}\neq X$ (i.e.
$B_{\overline{\gamma}_i},C_{\gamma_j}$ are just overlapping and do
not cover the whole set $X$). Then, the lower probability measure
of the credal set $\P_{\delta, \pi}$ is not $2-$monotone.
\end{prop}
The proof of Proposition~\ref{prop:CloudBel} can be found in the
appendix. It comes down to showing that for any non-comonotonic
cloud with a pair $B_{\overline{\gamma}_i},C_{\gamma_j}$ of events
satisfying the proposition, the situation exhibited in the above
example always occurs, namely the existence of two subsets of the
form $B_{\overline{\gamma}_i}$ and $C_{\gamma_j}^c$ for which
2-monotonicity fails. This indicates that random sets do not
generalize such non-comonotonic clouds.  It suggests that such non-comonotonic clouds are
likely to be less tractable when processing uncertainty: for
instance, simulation of such clouds via sampling methods  will be difficult to implement, and the computation of
lower/upper expectation too (since Choquet integral cannot be computed from credal sets  for
measures failing 2-monotonicity).

Note that comonotonic clouds
and clouds described by Proposition~\ref{prop:CloudBel} cover a
large number of possible discrete clouds, but that there remains
some "small" subfamilies, i.e. those non-comonotonic clouds for
which $\forall i, j, B_{\overline{\gamma}_i} \cap C_{\gamma_j} \in
\{B_{\overline{\gamma}_i},C_{\gamma_j},\emptyset \}$, or
$B_{\overline{\gamma}_i} \cup C_{\gamma_j}=X$. As such families
are very peculiar, we do not consider them further here.

\subsection{Outer approximation of a non-comonotonic cloud}
\label{sec:NoncoCloudOuter}

We provide, in this section and the next one, some practical means
to compute guaranteed outer and inner approximations of the exact
probability bounds induced by a non-comonotonic cloud, eventually
leading to an easier handling of such clouds.

Given a cloud $[\delta,\pi]$, we have proven that $\P_{[\delta,\pi]}
= \P_{\pi} \cap \P_{1-\delta}$, where $\pi$ and $1-\delta$ are
possibility distributions. As a consequence, the upper and lower
probabilities of $\P_{[\delta,\pi]}$ on any event can be bounded
from above (resp. from below), using the possibility measures and
the necessity measures induced by $\pi$ and $\overline{\pi}  =
1-\delta$. 
The following bounds, originally considered by
Neumaier~\cite{Neumaier04}, provide, for all event $A$ of $X$, an
outer approximation of the range of $P(A)$:
\begin{equation}
\label{eq:couterappr} \max(N_{\pi}(A), N_{1-\delta}(A))\!\! \leq
\underline{P}(A) \leq \! P(A)\! \leq \overline{P}(A) \leq
\!\! \min(\Pi_{\pi}(A), \Pi_{1-\delta}(A)),
\end{equation}
where $\underline{P}(A),\overline{P}(A)$ are the lower and upper
probabilities induced by $\P_{[\delta,\pi]}$. Remember that
probability bounds generated by possibility distributions alone are
of the form $[0,\beta]$ or $[\alpha,1]$. Using a cloud and applying
Equation (\ref{eq:couterappr}) lead to tighter bounds of the form
$[\alpha,\beta] \subset [0,1]$, while remaining simple to compute.
Nevertheless, these bounds are not, in general, the tightest ones enclosing $P(A)$ induced by $\P_{[\delta,\pi]}$, as the next
example shows:
\begin{exmp}
\label{exmp:Outerbounds} Let $[\delta,\pi]$ be a cloud defined on a
set $X$, such that distributions $\delta$ and $\pi$ takes up to
four different values on elements $x$ of $X$ (including 0 and 1).
These values are such that $0=\gamma_0 < \gamma_1 < \gamma_2 <
\gamma_3=1$, and the distributions $\delta,\pi$ are such that
\begin{displaymath}
\begin{array}{ccl@{\hspace{4mm}}ccl}
\pi(x) & = & 1 \mbox{~if~} x \in B_{\overline{\gamma}_{2}};  & \delta(x) & = & \gamma_2 \mbox{~if~} x \in C_{\gamma_2};\\
           & = &  \gamma_2 \mbox{~if~} x \in B_{\overline{\gamma}_{1}} \setminus B_{\overline{\gamma}_{2}};   &      & = & \gamma_1 \mbox{~if~} x \in C_{\gamma_1} \setminus C_{\gamma_2};\\
           & = &  \gamma_1 \mbox{~if~} x \not\in B_{\overline{\gamma}_{1}}.     &      & = & 0 \mbox{~if~} x \not\in C_{\gamma_1}.
\end{array}
\end{displaymath}
Since $P(B_{\overline{\gamma}_{1}}) \geq 1- \gamma_1$ and
$P(C_{\gamma_2}) \leq 1- \gamma_2$, from Equations
(\ref{eq:cloudcons}), we can check that
$\underline{P}(B_{\overline{\gamma}_{1}} \setminus C_{\gamma_2}) =
\underline{P}(B_{\overline{\gamma}_{1}} \cap C^c_{\gamma_2}) =
\gamma_2 - \gamma_1$. Now, by definition of a necessity measure,
$N_{\pi}(B_{\overline{\gamma_{1}}} \cap C^c_{\gamma_2}) =
\min(N_{\pi}(B_{\overline{\gamma}_{1}}),N_{\pi}(C^c_{\gamma_2})) =
0$ since $\Pi_{\pi}(C_{\gamma_2})= 1$ because $C_{\gamma_2}
\subseteq B_{\overline{\gamma}_{1}}$ and
$\Pi_{\pi}(B_{\overline{\gamma}_{1}})=1$. Considering distribution
$\delta$, we can have $N_{1-\delta}(B_{\overline{\gamma}_{1}} \cap
C^c_{\gamma_2}) =
\min(N_{1-\delta}(B_{\overline{\gamma}_{1}}),N_{1-\delta}(C^c_{\gamma_2}))
= 0$ since $N_{1-\delta}(B_{\overline{\gamma}_{1}}) =
\Delta_{\delta}(B_{\overline{\gamma}_{1}}^c)  = 0$ and $C_{\gamma_1}
\subseteq B_{\overline{\gamma}_{1}}$ (which means that the elements
$x$ of $X$ that are in $B_{\overline{\gamma}_{1}}^c$ are such that
$\delta(x)=0$). Equation (\ref{eq:couterappr}) can thus result in a
trivial lower bound (i.e. equal to 0), different from
$\underline{P}(B_{\overline{\gamma}_{1}} \cap C^c_{\gamma_2})$.
\end{exmp}
Bounds given by Equation (\ref{eq:couterappr}), are the main
motivation for clouds, after Neumaier~\cite{Neumaier04}. Since these
bounds are, in general, not optimal, Neumaier's claim that they are only vaguely
related to Walley's previsions or to random sets is not surprising.
Equation (\ref{eq:couterappr}) appears less useful in the case of
comonotonic clouds, for which optimal lower and upper probabilities of
events can be more easily computed (see Remark 3.7 in \cite{DesterckeAll07IJAR}).

\subsection{Inner approximation of a non-comonotonic cloud}
\label{sec:NoncoCloudInner}

The previous outer approximation is easy to compute and allows to
clarify some of Neumaier's claims. Nevertheless, it is still
unclear how to practically use these outer bounds in subsequent
treatments (e.g., propagation, fusion). The inner approximation of a
cloud $[\delta,\pi]$ proposed now is a random set, which is easy to
exploit in practice. This inner approximation is obtained as follows:
\begin{prop}
\label{prop:cinnerapp} Let $[\delta,\pi]$ be a non-comonotonic cloud
defined on $X$. Let us then define, for $j=1,\ldots,M$, the
following random set:
\begin{displaymath} \left\{ \begin{array}{c} E_j=\{ x \in X
| (\pi(x) \geq \gamma_j)
\wedge (\delta(x) < \gamma_j) \} \\
m(E_j)=\gamma_j - \gamma_{j-1}
\end{array} \right.
\end{displaymath}
where $0=\gamma_0 < \ldots \gamma_j < \ldots < \gamma_M=1$ are the
distinct values taken by $\delta,\pi$ on elements of $X$, $E_j$ are
the focal elements with masses $m(E_j)$ of the random set. This
random set is an inner approximation of $[\delta,\pi]$, in the sense
its credal set $\P_{Bel}$ is included in $ \P_{[\delta,\pi]}$.
\end{prop}
In the case of non-comonotonic clouds satisfying
Proposition~\ref{prop:CloudBel}, the inclusion is strict. This inner
approximation appears to be a natural candidate, since on events of
the type  $$\{
B_{\overline{\gamma_i}},C_{\gamma_i},B_{\overline{\gamma_i}}\setminus
C_{\gamma_j}| i=0,\ldots,M;j=0,\ldots,M;i \leq j \} $$, it gives
optimal bounds, and it is exact when the cloud $[\delta,\pi]$ is
comonotonic.

\section{Clouds and probability intervals}
\label{sec:probintclouds}

There is no direct relationship between clouds and probability
intervals~\cite{CamposAll94}. Nevertheless, we can study how to transform a set of
probability intervals into a cloud. Such transformations can be
useful when one wishes to work with clouds but information is given
in terms of probability intervals. There are mainly two paths that
can be followed to do this transformation:
\begin{itemize}
\item the first uses the representability of clouds by pairs
of possibility distributions, and extends existing transformations
of probability intervals into a single possibility distribution.
\item The second uses the equivalence between generalized
p-boxes and comonotonic clouds
\end{itemize}

\subsection{Exploiting probability-possibility transformations}
\label{sec:ProbintCloudClassical}

The problem of transforming a probability distribution into a
quantitative possibility distribution has been addressed by many
authors~\cite{DuboisAll00c}. A consistency
principle between (precise) probabilities and possibility
distributions was first informally stated by Zadeh~\cite{Zadeh78}:
what is probable should be possible. It was later translated by
Dubois and Prade~\cite{DuboisPrade80,DuboisAll91} as the following mathematical
constraint. Given a possibility distribution $\pi$ obtained by the
transformation of a probability measure $P$, this distribution
should be such that, for all events $A$ of $X$, we have $P(A) \leq
\Pi(A)$, with $\Pi$ the possibility measure of $\pi$ which is said
to dominate $P$. There are multiple possibility distributions
satisfying this requirement, and Dubois and
Prade~\cite{DuboisAll04,DuboisAll91} proposed to add the following
ordinal equivalence constraint, such that for two elements $x,y$ in
$X$
\begin{displaymath}
p(x) \leq p(y) \iff \pi(x) \leq \pi(y)
\end{displaymath}
and to choose the least specific possibility distribution ($\pi'$ is
more specific than $\pi$ if $\pi'\leq\pi$) respecting these two
constraints.

The unique solution~\cite{DuboisPrade80} is as follows: let us
consider probability masses such that $ p_1 \leq \ldots \leq
p_n $ with $p_j=p(x_j)$. When all probabilities are different, Dubois and
Prade probability-possibility transformation can be formulated as
\begin{displaymath}
\pi_i = \sum_{j=1}^i p_j
\end{displaymath}
with $\pi_i=\pi(x_i)$.  When some elements have equal probability,
the above equation must be used on the ordered partition induced by
the probability weights, using uniform probabilities inside each
element of the partition.

Reversing the ordering of the $p_i$'s in the above formula yields
another possibility distribution $\overline{\pi}_i = \sum_{j=i}^n
p_j$, with $\overline{\pi}_i=\overline{\pi}(x_i)$. Letting
$\delta=1-\overline{\pi}$, distribution $\delta$ is of the form
$\delta_\pi$ introduced in section \ref{Nonemptycredal},  that is,
$[\delta,\pi]$ is a cloud such that $\delta_i=\pi_{i-1}$ for all
$i>1$, with $\delta_1=0$ and $\delta_i=\delta(x_i)$. It is precisely the
tightest cloud containing $P$, in the sense that $\P_\pi \cap
\P_{\overline{\pi}} = \{P\}$. This shows that, at least when
probability masses are precise, transformation into possibility
distributions can be extended to get a second possibility
distribution such that this pair of distributions is equivalent to a
cloud that singles out $P$ exactly.

When working with imprecise probability assignments, i.e. with a probability interval $L$,\footnote{A probability interval on a space
$X$ is a tuple of intervals  $\{[l(x), u(x)] l x \in X\}$ enclosing the probabilities $p(x), x \in X$.}
 a partial order  $ \leq_L$ (actually, an interval
order) is induced by probability
weights on $X$ and defined by:
\begin{displaymath}
x \leq_L y \iff u(x) \leq l(y)
\end{displaymath}
and two elements $x,y$ are incomparable if intervals
$[l(x),u(x)],[l(y),u(y)]$ intersect. The problem of transforming a
probability interval into an outer-approximating
possibility distribution by extending Dubois and Prade
transformation is studied in detail by Masson and
Denoeux~\cite{MassonDenoeux06}. We first recall their method, before
proposing its extension to clouds.

Let
$\mathcal{C}_L$ be a set of linear extensions of the partial order
$\leq_L$: a linear extension $<_l \in \mathcal{C}_L$ is a linear ordering
of $X$ compatible with the partial order $\leq_L$. Let $\sigma_l$ be
the permutation such that $\sigma_l(x)$ is the rank of element $x$
in the linear extension $<_l$. Given the partial order  $ \leq_L$, Masson and
Denoeux~\cite{MassonDenoeux06} propose the following procedure:
\begin{enumerate}
\item For each linear order $<_l \in \mathcal{C}_L$ and each element $x$,
solve
\begin{equation}
\label{eq:probintposs1} \pi^l(x) = \max\limits_{\{p(y) | y \in X \}}
\sum\limits_{\sigma_l(y) \leq \sigma_l(x)} p(y)
\end{equation}
under the constraints
\begin{displaymath}
\left\{ \begin{array}{c}
  \sum\limits_{x \in X} p(x) = 1 \\
  \forall x \in X, \; l(x) \leq p(x) \leq u(x) \\
  p(\sigma_l^{-1}(1)) \leq p(\sigma_l^{-1}(2)) \leq \ldots \leq
p(\sigma_l^{-1}(n))
\end{array} \right.
\end{displaymath}
\item The most informative distribution $\pi$ dominating all distributions
$\pi^l$ is:
\begin{equation}
\label{eq:probintposs} \pi(x) = \max_{<_l \in \mathcal{C}} \pi^l(x).
\end{equation}
\end{enumerate}
This procedure ensures that the resulting possibility distribution
$\pi$ outer-approximate $\P_L$ (i.e. $\P_{L} \subseteq \P_{\pi}$).

Now, consider that the possibility distribution $\pi$ given by
Equation (\ref{eq:probintposs}) is the upper distribution of a cloud
$[\delta,\pi]$. To extend above procedure, we have to build a second
possibility distribution $\pi_{\delta}$ such that $\P_L \subseteq
\P_{\pi_{\delta}}$ and such that the pair $[1-\pi_{\delta},\pi]$
defines a cloud. To achieve this, we propose to use the same method
as Masson and Denoeux~\cite{MassonDenoeux06}, simply reversing the
inequality under the summation sign in Equation
(\ref{eq:probintposs1}). The procedure that builds $\pi_{\delta}$ then
becomes
\begin{enumerate}
\item For each order $<_l \in \mathcal{C}_L$ and each element $x$,
solve
\begin{align} \label{lowerdis}
\pi_{\delta}^l(x) & = \max\limits_{\{p(y) | y \in X
\}} \sum\limits_{\sigma_l(x) \leq \sigma_l(y)} p(y) \\
& = 1 - \min\limits_{\{p(y) | y \in X \}} \sum\limits_{\sigma_l(y) <
\sigma_l(x)} p(y) = 1 - \delta^l(x)
\end{align}
with the same constraints as in the first transformation.
\item The most informative distribution dominating all distributions
$\pi_{\delta}^l(x)$ is:
\begin{equation}
\label{eq:probintdelta} \pi_{\delta}(x) = 1 - \delta(x) = \max_{<_l
\in \mathcal{C}} \pi_{{\delta}}^l(x).
\end{equation}
\end{enumerate}
And we can check the following property:
\begin{prop}
\label{prop:ProbIntClouds} Given probability interval $L$,
the cloud $[1-\pi_{\delta},\pi]$ built from the two possibility
distributions $\pi_{\delta},\pi$ obtained via the above procedures
is such that the induced credal set $\P_{[1-\pi_{\delta},\pi]}$
outer-approximate $\P_L$. In the degenerate case of a precise
probability distribution, this cloud contains this distribution
only.
\end{prop}
\begin{proof}
The two possibility distributions $\pi,\pi_{\delta}$ are such that
$\P_L \subset \P_{\pi} $ and $\P_L \subset \P_{\pi_{\delta}}$ by
construction, so $\P_L \subset (\P_{\pi} \cap \P_{\pi_{\delta}})$.
The final result is thus more precise than a single possibility
distribution dominating $\P_L$. When $L$ reduces to a precise probability distribution
$\{p\}$, the transformations give the following possibility
distributions (elements of $X$ are ordered in accordance with the
order of probability masses):
\begin{displaymath}
\pi(x_i) = \sum_{j \leq i} p_j
\end{displaymath}
and
\begin{displaymath}
\pi_{\delta}(x_i) = \sum_{j \geq i} p_j = 1 - \sum_{j < i} p_j = 1 -
\delta(x_i) = 1 - \pi(x_{i-1}).
\end{displaymath}
Hence, the only probability distribution in the cloud $[\delta,
\pi]$ is given by $p_i = \pi(x_i) - \pi(x_{i-1})$.
\end{proof}

So, this method
constructs a cloud outer-approximating any
probability interval. It directly extends known methods used in
possibility theory.

\begin{exmp}
\label{exmp:IntCloud} Let us take the same probability
interval as in the example given by Masson and
Denoeux~\cite{MassonDenoeux06}, on the set $X=\{w,x,y,z\}$, and
summarized in the following table \begin{center}
\begin{tabular}{c@{\hspace{10mm}}cccc}
   &  $w$ & $x$ & $y$ & $z$ \\ \hline
  $l$  & 0.10  & 0.34 & 0.25 & 0  \\
   $u$  & 0.28 & 0.56 &  0.46 & 0.08  \\ \hline
\end{tabular}
\end{center}
The partial order is given by $L_y < L_x ; L_z < \{L_x,L_w,L_y\} $.
There are three possible linear extensions $<_l \in \mathcal{C}_L$
\begin{displaymath}
<_l^1  =  (L_z,L_w,L_y,L_x); \quad <_l^2  =  (L_z,L_w,L_x,L_y);
\quad <_l^3  =  (L_z,L_y,L_w,L_x)
\end{displaymath}
corresponding to the following $\pi_{{\delta}}$'s:
\begin{center}
\begin{tabular}{c@{\hspace{10mm}}c@{\hspace{5mm}}c@{\hspace{5mm}}c@{\hspace{5mm}}c}
 $<_l^i$ & $\pi_{\delta}(w)$ & $\pi_{\delta}(x)$ &
$\pi_{\delta}(y)$ & $\pi_{\delta}(z)$ \\ \hline
 1 & 1 & 0.16 & 0.63 & 1 \\
 2 & 1 & 0.9 & 0.46 & 1 \\
 3 & 0.75 & 0.5 & 1 & 1 \\ \hline
$\max$ & 1 & 0.9 & 1 & 1 \\
\end{tabular}
\end{center}
and, finally, the obtained cloud is:
\begin{center}
\begin{tabular}{c@{\hspace{10mm}}c@{\hspace{5mm}}c@{\hspace{5mm}}c@{\hspace{5mm}}c}
 & $w$ & $x$ & $y$ & $z$ \\ \hline
$\pi$  & 0.64 & 1 & 1 & 0.08 \\
 $\delta$ & 0 & 0.1 & 0 & 0 \\ \hline
\end{tabular}
\end{center}
where $\pi$ is the possibility distribution obtained by Masson and
Denoeux~\cite{MassonDenoeux06} using their method. Note that the
cloud is only a little more informative than the upper distribution
taken alone (indeed, the only added constraint is that $p(x) \leq
0.9$).
\end{exmp}

\subsection{Using generalized p-boxes}
\label{sec:ProbintCloudGenpbox}

Since generalized p-boxes and comonotonic clouds are equivalent
representations, we can directly use transformations of probability
intervals into generalized p-boxes (using Equations (14) in \cite{DesterckeAll07IJAR}) to get an outer-approximating comonotonic cloud.
Consider the following example:
\begin{exmp}
\label{exmp:IntCloud2} Let us consider the same probability
intervals as in example~\ref{exmp:IntCloud} and the following order
relationship $R$ on the elements: $z <_R w <_R y <_R x$. The
comonotonic cloud equivalent to the generalized p-box associated to
this order is:
\begin{center}
\begin{tabular}{c@{\hspace{10mm}}c@{\hspace{5mm}}c@{\hspace{5mm}}c@{\hspace{5mm}}c}
  & $w$ & $x$ & $y$ & $z$ \\ \hline
$\overline{F} = \pi$ & 0.36 & 1 & 0.66 & 0.08 \\
 $\underline{F}$ & 0.1 & 1 & 0.44 & 0 \\
 $\delta$ & 0 & 0.44 & 0.1 & 0 \\ \hline
\end{tabular}
\end{center}
\end{exmp}
And, using related results in the  companion paper \cite{DesterckeAll07IJAR}, we know that the
credal set $\P_{[\delta,\pi]}$ induced by this cloud is such that
$\P_L \subseteq \P_{[\delta,\pi]}$ and that we can recover the
information modeled by a probability interval by means
of at most $|X|/2$ clouds built by this method (Proposition 3.8 in \cite{DesterckeAll07IJAR}).


Both proposed methods transform a probability interval $L$ 
into a cloud $[\delta,\pi]$ outer-approximating $L$ (in the sense
that $\P_L \subset \P_{[\delta,\pi]}$, and in the case of a precise
probability distribution, each method recovers it exactly.

However, if we compare the clouds resulting from
Examples~\ref{exmp:IntCloud} and~\ref{exmp:IntCloud2}, it is clear
that the second method
(Example~\ref{exmp:IntCloud2}) is more precise than the first one (Example~\ref{exmp:IntCloud}).
Moreover,
using the first method, it is in general impossible to recover the
information provided by the original probability
interval. This shows that the first method can be very
conservative. This is mainly due to the fact that even if it
considers every possible ordering of elements, it is only based on
the partial order induced by the probability interval. Thus, if a
natural ordering of elements exists, the second method seems to
be preferable. Otherwise, it is harder to justify the fact of
considering one particular order rather than another one, and the
first method should be applied. In this case, one has to be aware
that a lot of information can be lost in the process. One may also
find out the ordering inducing the most precise comonotonic cloud,
but this question remains open.

\section{Continuous clouds on the real line}
\label{sec:Cloudrealline}

It often happens that uncertainty is defined on the real
line. It is thus important to know if results obtained so far can be
extended to continuous settings. In the following, we consider clouds defined on a
bounded interval $[\underline{r}, \overline{r}]$.

First, let us recall that, as in the discrete case, a cloud
$[\delta,\pi]$ defined on the real line is a pair of distributions
such that, for any element $r \in \mathbb{R}$, $[\delta(r),\pi(r)]$
is an interval and there is an element $r$ for which $\delta(r)=0$,
and another $r'$ for which $\pi(r')=1$. \emph{Thin clouds}
($\pi=\delta$) and \emph{fuzzy clouds} ($\delta=0$) have the same
definition as in the case of finite set. The credal set
$\P_{[\delta,\pi]}$ induced by a cloud on the real line is such
that:
\begin{equation}
\label{contcloud} \P_{[\delta,\pi]}=\{P | P( \{r \in \mathbb{R}
,\delta(r) \geq \alpha \} ) \leq 1- \alpha \leq P( \{r \in
\mathbb{R} ,\pi(r) > \alpha \} ) \},
\end{equation}
where $P$ is a $\sigma$-additive probability
distribution\footnote{To avoid mathematical subtleties that would
require special care, we restrict ourselves to
$\sigma$-additive probability distributions rather than considering
finitely additive ones.}.
\subsection{General results}
As Proposition 2.5 in \cite{DesterckeAll07IJAR} has been proven for very
general spaces~\cite{CousoAll01}, results whose proof is based on
this proposition directly extend to models on the real line.
Similarly, the proof of Proposition~\ref{lem:CloudBel1} extends
directly to continuous models on the real line.
Hence, the following statements still hold:
\begin{itemize}
\item if $[\delta,\pi]$ is a cloud, $1-\delta,\pi$ are possibility
distributions, and $\P_{[\delta,\pi]}=\P_{1-\delta} \cap \P_{\pi}$,
\item if $[\underline{F},\overline{F}]$ is a generalized p-box
defined on the reals, then
$\P_{[\underline{F},\overline{F}]}=\P_{\pi_{\underline{F}}} \cap
\P_{\pi_{\overline{F}}}$ with, for all $r \in \mathbb{R}$:
$$\pi_{\overline{F}}(r)=\overline{F}(r)$$ and
$$\pi_{\underline{F}}(r)= 1 - \sup\{\underline{F}(r') | r' \in \mathbb{R};\underline{F}(r') < \underline{F}(r) \}$$ with $\pi_{\underline{F}}(\underline{r})=0$.
\item {\em generalized} p-boxes and comonotonic clouds are equivalent representation
\end{itemize}

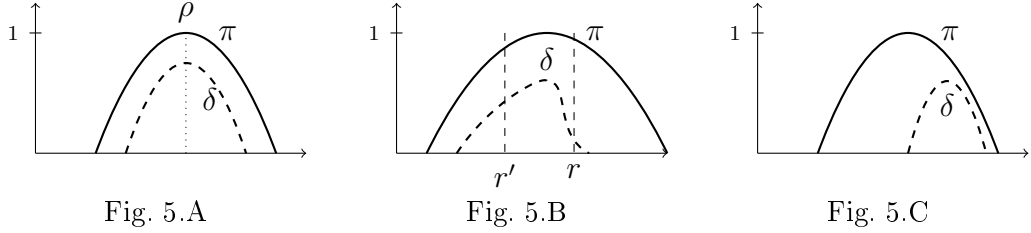
\begin{figure}
\begin{center}
\begin{tikzpicture}

\begin{scope}[scale=0.8]
\draw[thick] (10mm,0mm) parabola[parabola height=20mm] (40mm,0mm);

\draw[dashed,thick] (15mm,0mm) parabola[parabola height=15mm]
 (35mm,0mm) ;
\draw (32mm,20mm) node {$\pi$} ; \draw (29mm,9mm) node{$\delta$} ;

\draw[->] (0,0) -- (45mm,0mm) ; \draw (-1mm,20mm) node[left]
{$\scriptstyle 1$} -- (1mm,20mm) ; \draw[->] (0,0) -- (0mm,25mm) ;

\draw[dotted] (25mm,0mm) -- (25mm,20mm) node[above] {$\rho$} ;

\draw (20mm,-10mm) node {\small Fig.~\ref{fig:cloudsonreals}.A};
\end{scope}

\begin{scope}[scale=0.8,xshift=120mm]
\draw[thick] (10mm,0mm) parabola[parabola height=20mm] (40mm,0mm)  ;
\draw[dashed,thick] (25mm,0mm) parabola[parabola height=12mm]
(38mm,0mm);
\draw (32mm,20mm) node {$\pi$} ; \draw (31.5mm,11.7mm) node[below]
{$\delta$} ;

\draw[->] (0,0) -- (45mm,0mm) ; \draw (-1mm,20mm) node[left]
{$\scriptstyle 1$} -- (1mm,20mm) ; \draw[->] (0,0) -- (0mm,25mm) ;


\draw (20mm,-10mm) node {\small Fig.~\ref{fig:cloudsonreals}.C};
\end{scope}

\begin{scope}[scale=0.8,xshift=60mm]
\draw[thick] (5mm,0mm) parabola[parabola height=20mm] (45mm,0mm)  ;
\draw[dashed,thick] (10mm,0mm) .. controls (16mm,9mm) and
(24.5mm,13mm) .. (25mm,12mm) .. controls (29mm,12mm) and (26mm,3mm)
.. (32mm,0mm) ;
\draw (33mm,20mm) node {$\pi$} ; \draw (25mm,12mm) node[above]
{$\delta$} ;

\draw[dashed] (29.5mm,20mm) -- (29.5mm,0mm) node[below] {$r$} ;
\draw[dashed] (18mm,20mm) -- (18mm,0mm) node[below] {$r'$} ;

\draw[->] (0,0) -- (45mm,0mm) ; \draw (-1mm,20mm) node[left]
{$\scriptstyle 1$} -- (1mm,20mm) ; \draw[->] (0,0) -- (0mm,25mm) ;


\draw (20mm,-10mm) node {\small Fig.~\ref{fig:cloudsonreals}.B};
\end{scope}



\end{tikzpicture}
\caption{Illustration of comonotonic (A), weakly comonotonic (B) and
non-comonotonic clouds (C) on the real line. }
\label{fig:cloudsonreals}
\end{center}
\end{figure}

Note that, for clouds on the real line, we can define a weaker
notion of comonotonicity: a (continuous) cloud $[\delta,\pi]$ is
said to be \emph{weakly} comonotonic if the sign of the derivative
of distributions $\delta,\pi$ is the same in every point $r$ of the
real line $\mathbb{R}$. Being \emph{weakly} comonotonic is not
sufficient to be equivalent to a generalized p-box, since if $\pi$
and $\delta$ are only weakly comonotonic, then it is possible to
find two values $r$ and $r'$ such that $\delta(r)<\delta(r')$ and
$\pi(r)>\pi(r')$. In this case, the (pre-)ordering jointly induced
by the two distributions is not complete, and the definition of
comonotonicity is not satisfied.
 Figures~\ref{fig:cloudsonreals}.A,~\ref{fig:cloudsonreals}.B and~\ref{fig:cloudsonreals}.C respectively illustrate the notion of
comonotonic, non-comonotonic and weakly comonotonic clouds on the
reals. Figure~\ref{fig:cloudsonreals}.A illustrates a comonotonic
cloud (and, consequently, a generalized p-box) for which elements
are ordered according to their distance to the mode $\rho$ (i.e., for
this particular cloud, two values $x,y$ in $\mathbb{R}$ are such
that $x <_{[\underline{F},\overline{F}]} y$ if and only if $|\rho-x| >
|\rho-y|$). Note that Figure~\ref{fig:cloudsonreals}.A is a good illustration of the potential use of a generalized p-box, as already noticed (see beginning of section 3 in companion paper~\cite{DesterckeAll07IJAR}).

We can now extend the propositions linking clouds and generalized
p-boxes with random sets. In particular, the following result
extends Proposition~\ref{prop:CloudBel} to the continuous case:
\begin{prop}
\label{prop:contcloud} Let the distributions $[\delta,\pi]$ describe
a continuous cloud on the reals and $\P_{[\delta,\pi]}$ be the
induced credal set. Then, the random set defined by the Lebesgue
measure on the unit interval $\alpha \in [0,1]$ and the multimapping
$\alpha \longrightarrow E_{\alpha}$ such that
\begin{displaymath}
E_{\alpha}=\{ r \in \mathbb{R} | (\pi(r) \geq \alpha) \wedge
(\delta(r) < \alpha) \}
\end{displaymath}
defines a credal set $\P_{Bel}$ inner-approximating
$\P_{\pi,\delta}$ ($\P_{Bel} \subseteq \P_{\pi,\delta}$).
\end{prop}
The proof can be found in the appendix. It comes down to using
sequences of discrete clouds outer- and inner-approximating
$[\delta,\pi]$ and converging to it, and then to consider
inner-approximations of those discrete clouds given by
Proposition~\ref{prop:cinnerapp}. This proposition has two
corollaries:
\begin{cor}
\label{cor:contRScomocloud} Let $[\delta,\pi]$ be a comonotonic
cloud with continuous distributions on the real line. Then the
credal set $\P_{[\delta,\pi]}$ is  also the credal set of a
continuous random set with uniform mass density, whose focal sets
are of the form, for $\alpha \in [0,1]$:
\begin{displaymath} E_{\alpha}=\{ r \in \mathbb{R} | (\pi(r) \geq \alpha) \wedge
(\delta(r) < \alpha) \}.
\end{displaymath}
\end{cor}
To obtain the result, simply observe that the inner-approximation of
Proposition~\ref{prop:cinnerapp} becomes exact for discrete comonotonic
clouds, which are special cases of random sets. In particular, this is true for the
sequences of discrete comonotonic clouds outer- and inner-approximating
$[\delta,\pi]$ and converging to it. So, this sequence of random sets converge to a continuous random set at the limit. Another interesting particular case is the one of uniformly
continuous p-boxes.
\begin{cor}
The credal set $\P_{[\underline{F},\overline{F}]}$ described by two
continuous and strictly increasing cumulative distributions
$\underline{F},\overline{F}$ forming a classical p-box on the reals
is equivalent to the credal set described by the continuous random
set with uniform mass density, whose focal sets are sets of the form
$[x(\alpha),y(\alpha)]$ where $x(\alpha) =
\overline{F}^{-1}(\alpha)$ and $y(\alpha) =
\underline{F}^{-1}(\alpha).$
\end{cor}
This is because strictly increasing continuous p-boxes are special
cases of comonotonic clouds (or, equivalently, of generalized
p-boxes). To check that, in this case, $E_\alpha=[x(\alpha),y(\alpha)]$, it suffices to consider the possibility distributions $\pi_{\underline{F}},\pi_{\overline{F}}$ and to check that $\inf_r\{\pi_{\overline{F}}(r)\geq \alpha)\}=x(\alpha)$ and that $\sup_r\{1-\pi_{\underline{F}}(r)<\alpha\}=y(\alpha)$.The strict increasingness property can be relaxed to
intervals where the  cumulative functions are constant, provided one
consider pseudo-inverses when building the continuous random set.

These results are interesting, for they can make the computation of
lower and upper expectations over continuous generalized p-boxes
easier. Another interesting point is that the framework
developed by Smets~\cite{Smets05} concerning belief functions on the
reals can be applied to comonotonic clouds. Also note that above
results extend and give alternative proofs to other results given by
Alvarez~\cite{Alvarez06} concerning continuous p-boxes.

\subsection{Thin clouds}
The case of thin clouds,  for which $\pi=\delta$, is interesting. In this case, constraints
(\ref{eq:cloudcons}) defining the credal set $\P_{[\delta,\pi]}$
reduce to $ P( \pi(x) \geq \alpha ) = P( \pi(x) > \alpha )= 1-
\alpha $ for all $\alpha \in (0, 1)$. As noticed earlier, when $X$
is finite, thin clouds define empty credal sets, but is no longer
the case when it is defined on the real line, as the following
proposition shows:
\begin{prop}
\label{prop:contthinclouds} If $\pi$ is a continuous possibility
distribution on the real line, then the credal set  $\P_{[\pi,\pi]} = \P_{\pi} \cap
\P_{1-\pi}$ is not empty.
\end{prop}
\begin{proof}[\textbf{Proof of Proposition~\ref{prop:contthinclouds}}]  Let $F(x) = \Pi((- \infty, x])$, with $x \in \mathbb{R}$. $F$ is the distribution
function of a probability measure $P_{\pi}$ such that for all $
\alpha \in [0,1]$, $P_{\pi}(\{x \in \mathbb{R} | \pi(x)
>\alpha \} ) = 1- \alpha$, where the sets
$\{ x \in \mathbb{R} | \pi(x) >\alpha \}$ form a continuous nested
sequence (see \cite{DuboisAll04} p. 285). Such a probability lies in
$\P_{\pi}$. Moreover, $$P_{\pi}( \{ x \in \mathbb{R} | \pi(x)
>\alpha \}  )= P_{\pi}( \{ x \in \mathbb{R} | \pi(x)
\geq \alpha \} )$$ due to uniform continuity of $\pi$. We also have
\\ \mbox{$P_{\pi}( \{ x \in \mathbb{R} | \pi(x)
>\alpha \}  ) =  1- \Pi( \{ x \in \mathbb{R} | \pi(x)
\geq \alpha \}^c) = 1 - \Delta( \{ x \in \mathbb{R} | \pi(x) \geq
\alpha \}) $} again due to uniform continuity. Since
\\ \mbox{$1 - \Delta( \{ x \in \mathbb{R} | \pi(x) \geq \alpha \})=
\sup_{ x| \pi(x) \geq \alpha} 1 - \pi(x)$}, this means $P_{\pi} \in
\P(1-\pi)$.
\end{proof}

A \emph{thin} cloud is a particular case of comonotonic cloud. It
induces a complete pre-ordering on the reals. If this pre-order is
linear, it means that for any $\alpha \in [0,1]$, there is only one
value $r \in \mathbb{R}$ for which $\pi(r)=\alpha$, and that
$\P_{\pi} \cap \P_{1-\pi}$ contains only one probability measure. In
particular, if the order is the natural order of real numbers, this
\emph{thin} cloud reduces to an usual cumulative distribution. When
the pre-order has ties, it means that for some $\alpha \in [0,1]$,
there are several values in $r \in\mathbb{R}$ such that
$\pi(r)=\alpha$. Using Corollary~\ref{cor:contRScomocloud}, we can
model the credal set $P_{\pi} \in \P(1-\pi)$ by the random set with
uniform mass density, whose focal sets are of the form
\begin{displaymath}
E_{\alpha}=\{r \in \mathbb{R} | \pi(r)=\alpha \}
\end{displaymath}
In this case, we can check that $Bel(\{r \in \mathbb{R} | \pi(r)
\geq \alpha \})  = 1 - \alpha$, in accordance with Equation
(\ref{eq:cloudcons}).

Finally, consider the specific case of a thin cloud modeled by an
unimodal distribution $\pi$ (formally, a fuzzy interval). In this
case, each focal set associated to a value $\alpha$ is a doubleton
$\{x(\alpha), y(\alpha)\}$ where $\{x | \pi(x) \geq \alpha \} =
[x(\alpha), y(\alpha)]$. Noticeable probability distributions that
are inside the credal set induced by such a \emph{thin} cloud are
the cumulative distributions $F_+$ and $F_-$ such that for all
$\alpha$ in $[0,1]$ $F_+^{-1}(\alpha)=x(\alpha)$ and
$1-F_-^{-1}(\alpha)=y(\alpha)$ (they respectively correspond to a
mass density concentrated on values $x(\alpha)$ and $y(\alpha)$).
All probability measures with cumulative functions of the form
$\lambda\cdot F_+ + (1- \lambda)\cdot F_-$ also belong to the credal
set (for $\lambda = \frac{1}{2}$, this distribution is obtained by
evenly dividing mass density between elements $x(\alpha)$ and
$y(\alpha)$). Other distributions inside this set are considered by
Dubois \emph{et al.}~\cite{DuboisAll04}.

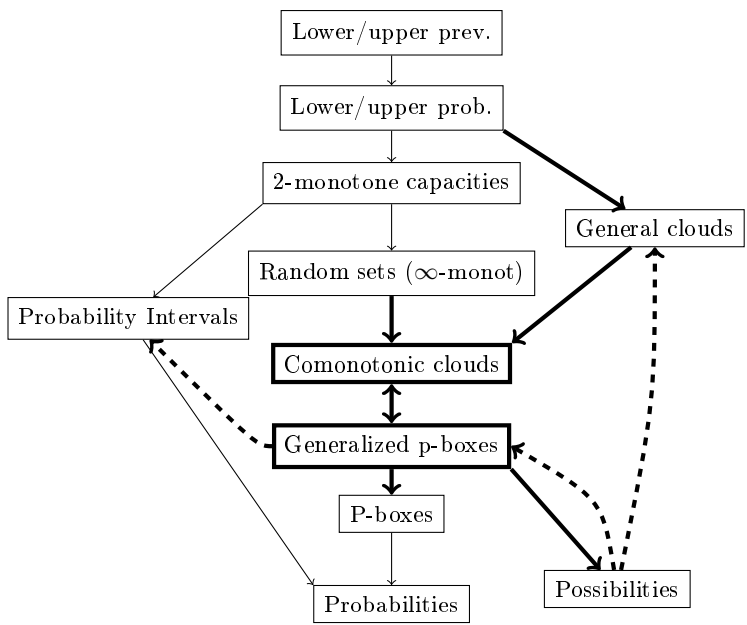
\begin{figure}
\begin{center}
\begin{tikzpicture}
 \node (ImpProb) at (0,-0.4) [draw] {\scriptsize{Lower/upper prev.}};
 \node (LowUpp) at (0,-1.4) [draw] {\scriptsize{Lower/upper prob.}};
 \node (2mon) at (0,-2.4) [draw] {\scriptsize{2-monotone capacities}} ;
 \node (randset) at (0,-3.6) [draw] {\scriptsize{Random sets ($\infty$-monot)}} ;
 \node (comclo) at (0,-4.8) [draw,ultra thick] {\scriptsize{Comonotonic clouds}} ;
 \node (GenPbox) at (0,-5.9) [draw,ultra thick] {\scriptsize{Generalized p-boxes}} ;
 \node (Pbox) at (0,-6.8) [draw] {\scriptsize{P-boxes}} ;
 \node (Prob) at (0,-8) [draw] {\scriptsize{Probabilities}} ;
 \node (ProbInt) at (-3.5,-4.2) [draw] {\scriptsize{Probability
 Intervals}} ;
 \node (NonCoClo) at (3.5,-3) [draw] {\scriptsize{General clouds}} ;
 \node (Poss) at (3,-7.8) [draw] {\scriptsize{Possibilities}} ;

 \draw[->,ultra thick] (LowUpp.south east) -- (NonCoClo) ;
 \draw[->,ultra thick] (NonCoClo) -- (comclo.north east) ;
 \draw[->] (ImpProb) -- (LowUpp)  ;
 \draw[->]  (LowUpp) -- (2mon) ;
 \draw[->]  (2mon) -- (randset)  ;
 \draw[->,ultra thick]  (randset) -- (comclo) ;
 \draw[<->,ultra thick] (comclo) -- (GenPbox) ;
 \draw[->,ultra thick] (GenPbox) -- (Pbox)  ;
 \draw[->]  (Pbox) -- (Prob) ;
 \draw[->,ultra thick] (GenPbox.south east) -- (Poss) ;
 \draw[->] (2mon.south west) -- (ProbInt)  ;
 \draw[->]  (ProbInt) -- (Prob.north west) ;
 \draw[dashed,->,ultra thick]  (Poss) .. controls
+(-0.2,1.2) .. (GenPbox.east) ; \draw[dashed,->,ultra thick]
(GenPbox.west) .. controls +(-0.2,0) .. (ProbInt) ;
\draw[dashed,->,ultra thick] (Poss) .. controls +(0.5,2.5) ..
(NonCoClo) ;

\end{tikzpicture}  \caption{Representation relationships: completed summary with clouds. $A \longrightarrow B$: B
  is a special case of A. $A \dashrightarrow B$: B is representable by A }\label{fig:summary}
\end{center}
\end{figure}

\section{Conclusion}


In this paper  Neumaier clouds are compared to other uncertainty
representations, including generalized p-boxes introduced in the
 companion paper \cite{DesterckeAll07IJAR}. Properties of the cloud formalism are explained in
the light of other representations. We are now ready to
complete Figure~\ref{fig:summarypre} with clouds. This completed
picture is given by Figure~\ref{fig:summary}. New relationships and
representations coming from this paper and its companion are in bold
lines.



The next step is to explore computational aspects of each formalism
as done by De Campos \emph{et al.}~\cite{CamposAll94} for
probability intervals. In particular, we need to answer the
following questions: how do we define operations of fusion,
marginalization, conditioning or propagation for each of these
models? Are the representations preserved after such operations, and
under which assumptions? What is the computational complexity of
these operations? Can the models presented here be easily elicited
or integrated? If many results already exist for random sets,
possibility distributions and probability intervals, few have been
derived for generalized p-boxes or clouds, due to their novelty. The
results presented in this paper and its companion can be helpful to
perform such a study. Recent applications of clouds to engineering
design problems~\cite{FuchsNeum08} indicate that this model can be
useful, and that such a study should be done to gain more insight
about the potential of such models. In particular, the mathematical
properties of comonotonic clouds appear to be quite attractive.
Our study thus indicates how clouds and generalized p-boxes can be
interpreted in the framework of other uncertainty theories. 


%

Another issue is to extend presented results to more general spaces,
to general lower/upper previsions or to cases not considered here
(e.g. continuous clouds with some discontinuities), possibly by
using existing results~\cite{Smets05,GDecoo05}. 

\section*{Appendix}

 We first recall a useful result
by Chateauneuf~\cite{Chateauneuf94} concerning the intersection of
credal sets induced by random sets. Consider two random sets $\{(F_i, q_{ i\cdot}) | i = 1, \dots k\}$ and $\{(G_j, q_{\cdot j})| j = 1, \dots l\}$ on $X$, with $F_i,G_j$ the focal elements, $q_{ i\cdot},q_{\cdot j}$ the corresponding masses and  $\mathcal{P}_F$ and $\mathcal{P}_G$ the induced credal sets. Consider then the set $\mathcal{Q}$ of all random sets $Q$ of the form $\{(F_i\cap G_j, q_{ ij})| i = 1, \dots k; , j = 1, \dots l\}$, with $F_i\cap G_j$  the focal sets and $q_{ij}$ the masses such that $q_{i\cdot} = \sum_{j=1}^l  q_{ij}  $ and  $q_{\cdot j} =  \sum_{i=1}^{k} q_{ij}$ with the constraint that $q_{ ij} = 0$ whenever $F_i\cap G_j = \emptyset$. Then the lower probability induced by the credal set $\mathcal{P}_F\cap \mathcal{P}_G$ is $$\underline{P}(A)=\min_{P \in \mathcal{P}_1\cap \mathcal{P}_2}P(A) = \min_{Q \in \mathcal{Q}}Bel_Q(A), \forall A \subseteq X,$$ where $Bel_Q$ is the belief function induced by the random set $Q$. 
\begin{proof}[\textbf{Proof of Proposition~\ref{prop:CloudBel}}]
We first state a short
Lemma allowing us to emphasize the idea behind the proof of the
latter proposition.
\begin{lem}
\label{lem:FourCloud} Let $(F_1,F_2),(G_1,G_2)$ be two pairs of sets
such that $F_1 \subset F_2$, \mbox{$G_1 \subset G_2$}, $G_1
\nsubseteq F_2$ and $G_1 \cap F_1 \neq \emptyset$. Let also
$\pi_F,\pi_G$ be two possibility distributions such that the
corresponding belief functions are defined by mass assignments
$m_F(F_1) = m_G(G_2) = \lambda$, $m_F(F_2) = m_G(G_1) = 1 -
\lambda$. Then, the lower probability of the non-empty credal set
$\P = \P_{\pi_F} \cap \P_{\pi_G}$ is not $2-monotone$.
\end{lem}

\begin{proof}[\textbf{Proof of Lemma~\ref{lem:FourCloud}}]
Chateauneuf's  result is applied to
the possibility distributions defined in Lemma~\ref{lem:FourCloud}.
The main idea is to exhibit two events and computing their lower probabilities, showing that $2$-monotonicity is
violated.
Consider the set $\mathcal{M}$ of matrices $M$ of the form
\begin{displaymath}
\begin{array}{c|cc}

& G_1 & G_2 \\
\hline
F_1 & m_{11} & m_{12} \\
 F_2 & m_{21} & m_{22} 
\end{array}
\end{displaymath}
 where
\begin{align*}
& m_{11} + m_{12} = m_{22} + m_{12} = \lambda \\
& m_{21} + m_{22} = m_{21} + m_{11} = 1 - \lambda  \\
& \sum m_{ij} = 1
\end{align*}
Each such matrix is a normalized (i.e. such that $m(\emptyset)=0$)
joint mass distribution for the random sets induced from possibility
distributions $\pi_F,\pi_G$, viewed as marginal belief functions.
Following Chateauneuf~\cite{Chateauneuf94}, for any event $E \subseteq X$, the lower probability
$\underline{P}(E)$ induced by the credal set $\P = \P_{\pi_F} \cap
\P_{\pi_G}$ is
\begin{equation}
\label{eq:fourcloudlowprob} \underline{P}(E) = \min_{M \in
\mathcal{M}} \sum_{(F_i \cap G_j) \subset E} m_{ij}
\end{equation}
Now consider the four events $F_1,G_1,F_1 \cap G_1, F_1 \cup G_1$.
Studying the relations between sets and the constraints on the
values $m_{ij}$, we can see that
\begin{displaymath}
 \underline{P}(F_1) = \lambda; \quad
 \underline{P}(G_1) = 1 - \lambda; \quad
 \underline{P}(F_1 \cap G_1) = 0.
\end{displaymath}
For $F_1 \cap G_1$, just consider the matrix $m_{12} = \lambda,
m_{21} = 1 - \lambda$. To show that the lower probability is not
even $2-$monotone, it is enough to show that
\mbox{$\underline{P}(F_1 \cup G_1) < 1 $}. To achieve this, consider
the following mass distribution
\begin{align*}
& m_{11} = \min(\lambda,1-\lambda) & m_{21} = 1 - \lambda - m_{11} & \\
& m_{12} = \lambda - m_{11} &  m_{22} = \min(\lambda,1-\lambda). &
\end{align*}
It can be checked that this matrix 
is in the set $\mathcal{M}$, and yields
\begin{align*}
P(F_1 \cup G_1) & =   m_{12} + m_{11} + m_{21} = m_{11} + \lambda - m_{11} + 1 - \lambda - m_{11}\\
& =  1 - m_{11} = 1 - \min(\lambda,1-\lambda) = \max(1 - \lambda,
\lambda) < 1
\end{align*}
since  $(F_2 \cap G_2) \nsubseteq (F_1 \cup G_1)$ (due to the fact
that $G_1 \nsubseteq F_2$). Then the inequality
$\underline{P}(F_1 \cup G_1) + \underline{P}(F_1 \cap G_1) <
\underline{P}(F_1) + \underline{P}(G_1)$  violates 2-monotonicity.
\end{proof}
To prove Proposition~\ref{prop:CloudBel}, we again use the result by
Chateauneuf~\cite{Chateauneuf94}, and we exhibit a pair of
events for which 2-monotonicity fails. Chateauneuf results are
 applicable to clouds, since possibility distributions are
equivalent to nested random sets. Consider a finite cloud described
by Equations (\ref{eq:cloudcons}) and the following matrix $Q$ of
weights $q_{ij}$
\begin{displaymath}
\begin{array}{c|ccccccc}
& C_{\gamma_1}^c  &  \cdots & C_{\gamma_j}^c & \cdot & C_{\gamma_{i+1}}^c & \cdots & C_{\gamma_M}^c  \\
\hline
B_{\overline{\gamma}_0} & q_{11}  & \ldots & q_{1j} & \cdot & q_{1(i+1)} & \ldots & q_{1M} \\
\vdots & \vdots  & \ddots & &  & \vdots&  & \vdots \\
B_{\overline{\gamma}_{j-1}} & q_{j1}  & \ldots & \mathbf{q_{jj}} & \cdot & \mathbf{q_{j(i+1)}} & \ldots & q_{jM} \\
\vdots & \vdots & \vdots & \vdots & & \vdots & \ddots & \vdots\\
B_{\overline{\gamma}_i} & q_{(i+1)1}  & \ldots & \mathbf{q_{(i+1)j}} & \cdot & \mathbf{q_{(i+1)(i+1)}} & \ldots & q_{(i+1)M} \\
\vdots & \vdots & \vdots & \vdots & & \vdots & \ddots & \vdots\\
B_{\overline{\gamma}_{M-1}} & q_{M1}  & \ldots & q_{Mj} & \cdot &
q_{M(i+1)} & \ldots & q_{MM}
\end{array}
\end{displaymath}
Respectively call $Bel_1$ and $Bel_2$ the belief functions
equivalent to the possibility distributions respectively generated
by
the collections of sets \\
\mbox{$\{B_{\overline{\gamma}_i}|i=0,\ldots,M-1 \}$} and
$\{C_{\gamma_i}^c|i=1,\ldots,M \}$. Using the fact that possibility
distributions can be mapped into random sets, we have
$m_1(B_{\overline{\gamma}_i})=\gamma_{i+1} - \gamma_{i}$ for
$i=0,\ldots,M-1$, and $m_2(C_{\gamma_j}^c)=\gamma_j - \gamma_{j-1}$
for $j=1,\ldots,M$. As in the proof of Lemma~\ref{lem:FourCloud}, we
consider every possible joint random set such that $m(\emptyset)=0$
built from the two marginal belief functions $Bel_1,Bel_2$.
Following Chateauneuf, let $\mathcal{Q}$ be the set of matrices $Q$
s.t.
\begin{align*}
& q_{i\cdot} = \sum_{j=1}^M  q_{ij}  =  \gamma_{i} - \gamma_{i-1}  \\
& q_{\cdot j} =  \sum_{i=1}^{M} q_{ij}  =  \gamma_{j} - \gamma_{j-1}  \\
& \textrm{ If } i,j  \textrm{ s.t. } B_{\overline{\gamma}_i} \cap
C_{\gamma_j}^c =   \emptyset \textrm{ then } q_{ij}  = 0
\end{align*}
and the lower probability of the credal set $\P_{[\delta,\pi]}$ on
event $E$ is such that
\begin{equation}
\label{eq:cloudlowprob} \underline{P}(E) = \min_{Q \in \mathcal{Q}}
 \sum_{(B_{\overline{\gamma}_i} \cap C_{\gamma_j}^c) \subset E} q_{ij}.
\end{equation}
Now, by hypothesis, there are at least two overlapping sets
$B_{\overline{\gamma}_i},C_{\gamma_j} \; i>j$ that are not included
in each other (i.e. $B_{\overline{\gamma}_i} \cap C_{\gamma_j}
\not\in \{ B_{\overline{\gamma}_i},C_{\gamma_j},\emptyset \}$). Let
us consider the four events
$B_{\overline{\gamma}_i},C_{\gamma_j}^c,B_{\overline{\gamma}_i} \cap
C_{\gamma_j}^c,B_{\overline{\gamma}_i} \cup C_{\gamma_j}^c$, which
are all different by hypothesis. Considering Equation
(\ref{eq:cloudlowprob}), the matrix $Q$ and the relations between
sets, inclusions $B_{\overline{\gamma}_m} \subset \ldots \subset
B_{\overline{\gamma}_0}$, $C_{\gamma_0}^c \subset \ldots \subset
C_{\gamma_m}^c$ and, for $i=0,\ldots,m$, $C_{\gamma_i} \subset
B_{\overline{\gamma}_i}$ imply:
\begin{displaymath}
\underline{P}(B_{\overline{\gamma}_i}) = 1-\gamma_i; \quad
\underline{P}(C_{\gamma_j}^c) = \gamma_{j}; \quad
\underline{P}(B_{\overline{\gamma}_i} \cap C_{\gamma_j}^c) = 0.
\end{displaymath}
For the last result, just consider the mass distribution $q_{kk} =
\gamma_{k-1} - \gamma_{k}$ for $k=1,\ldots,m$.

Next, consider event $B_{\overline{\gamma}_i} \cup C_{\gamma_j}^c$
(which is different from $X$ by hypothesis), and let them play the
role of events $F_1,G_1$ in Lemma~\ref{lem:FourCloud}. Suppose all
masses are such that $q_{kk} = \gamma_{k-1} - \gamma_{k}$, except
for masses (in boldface in the matrix) $q_{jj}, q_{(i+1)(i+1)}$.
Then, $C_{\gamma_j}^c \subset C_{\gamma_{i+1}}^c$,
$B_{\overline{\gamma}_i} \subset B_{\overline{\gamma}_{j-1}}$,
$C_{\gamma_j}^c \nsubseteq B_{\overline{\gamma}_{j-1}}$ by
definition of a cloud and $B_{\overline{\gamma}_i} \cap
C_{\gamma_j}^c \neq \emptyset$ by hypothesis. Finally, using
Lemma~\ref{lem:FourCloud}, consider the mass distribution
\begin{align*}
 & q_{(i+1)(i+1)} = \gamma_{i+1} - \gamma_{i} - q_{(i+1)j} & q_{(i+1)j} = \min(\gamma_{i+1} - \gamma_{i},\gamma_j -
\gamma_{j-1})  \\
& q_{j(i+1)} = \min(\gamma_{i+1} - \gamma_{i},\gamma_j -
\gamma_{j-1}) &  q_{jj} = \gamma_j -
\gamma_{j-1} - q_{(i+1)j}.
\end{align*}
It always gives a matrix in the set $\mathcal{Q}$. By considering
every subset of $B_{\overline{\gamma}_i} \cup C_{\gamma_j}^c$, we
thus get the following inequality
\begin{displaymath}
\underline{P}(B_{\overline{\gamma}_i} \cup C_{\gamma_j}^c) \leq
\gamma_{j-1} + 1 - \gamma_{i+1} + \max(\gamma_{i+1} -
\gamma_{i},\gamma_j - \gamma_{j-1}).
\end{displaymath}
And, similarly to what was found in Lemma~\ref{lem:FourCloud}, we
get
\begin{displaymath}
\underline{P}(B_{\overline{\gamma}_i} \cup C_{\gamma_j}^c) +
\underline{P}(B_{\overline{\gamma}_i} \cap C_{\gamma_j}^c) <
\underline{P}(B_{\overline{\gamma}_i}) +
\underline{P}(C_{\gamma_j}^c),
\end{displaymath}
which shows that the lower probability is not $2-$monotone.
\end{proof}

\begin{proof}[\textbf{Proof of Proposition~\ref{prop:cinnerapp}}]
First, we know that the random set given in
Proposition~\ref{prop:cinnerapp} is equivalent to
\begin{displaymath}
\left\{ \begin{array}{c} E_j=
B_{\overline{\gamma}_{j-1}} \setminus C_{\gamma_j} = B_{\gamma_j} \setminus C_{\gamma_j} \\
m(E_j)=\gamma_j - \gamma_{j-1}
\end{array} \right.
\end{displaymath}
Now, if we consider the matrix given in the proof of
Proposition~\ref{prop:CloudBel}, this random set comes down, for
$i=1,\ldots,M$ to assign masses $q_{ii}=\gamma_i-\gamma_{i-1}$.
Since this is a legal assignment, we are sure that for all events $E
\subseteq X$, the belief function of this random set is such that
$Bel(E) \geq \underline{P}(E)$, where $\underline{P}$ is the lower
probability induced by the cloud. The proof of
Proposition~\ref{prop:CloudBel} shows that this inclusion is strict
for clouds satisfying the latter proposition (since the lower
probability induced by such clouds is not $2$-monotone). \end{proof}

\begin{proof}[\textbf{Proof of Proposition~\ref{prop:contcloud}}]

We build outer and inner approximations of the continuous random set
that converge to the belief measure of the continuous random set,
while the corresponding clouds of which they are inner
approximations themselves converge to the uniformly continuous
 cloud.

First, consider a finite collection  \mbox{$0 = \alpha_0 < \alpha_1 < \ldots < \alpha_n = 1$} of equidistant levels $\alpha_i$
(\mbox{$\alpha_{i-1} - \alpha_{i} = 1/n, \forall i=1,\ldots,n$}).
Then, consider the following discrete non-comonotonic clouds
$[\underline{\delta}_n,\underline{\pi}_n]$,
$[\overline{\delta}_n,\overline{\pi}_n]$ that are respectively outer
and inner approximations of the cloud $[\delta,\pi]$: for every
value $r$ in $\mathbb{R}$, do the following transformation
\begin{displaymath}
\begin{array}{ccc}
\pi(r)=\alpha \textrm{ with } \alpha \in [\alpha_{i-1},\alpha_i]
\textrm{ then } \underline{\pi}_n(r)= \alpha_i  \textrm{ and }  \overline{\pi}_n(r)= \alpha_{i-1} \\
\delta(r)=\alpha' \textrm{ with } \alpha' \in
[\alpha_{j-1},\alpha_j]   \textrm{ then }  \underline{\delta}_n(r)=
\alpha_{j-1}  \textrm{ and }  \overline{\delta}_n(r)= \alpha_{j}
\end{array}
\end{displaymath}
This construction is illustrated in Figure~\ref{fig:noncocloudappr}
for the particular case when both $\pi$ and $\delta$ are unimodal.
In this particular case, for $i=1,\ldots,n$
\begin{align*}
 \{x \in \mathbb{R}|\underline{\pi}(x) \geq \alpha \} = [x(\alpha_{i-1}),
y(\alpha_{i-1})] & \textrm{ with } \alpha \in [\alpha_{i-1},\alpha_{i}]    \\
 \{x \in \mathbb{R}|\underline{\delta}(x) > \alpha \} = [u(\alpha_{i}), v(\alpha_{i})]
& \textrm{ with } \alpha \in [\alpha_{i-1},\alpha_{i}]
\end{align*}
\begin{align*}
\{x \in \mathbb{R}|\overline{\pi}(x) \geq \alpha \} =
[x(\alpha_{i}),
y(\alpha_{i})] & \alpha \in [\alpha_{i-1},\alpha_{i}]  \\
\{x \in \mathbb{R}|\overline{\delta}(x) > \alpha \} =
[u(\alpha_{i-1}), v(\alpha_{i-1})] & \alpha \in
[\alpha_{i-1},\alpha_{i}]
\end{align*}

\begin{figure}
\begin{center}
\includegraphics[width=8cm]{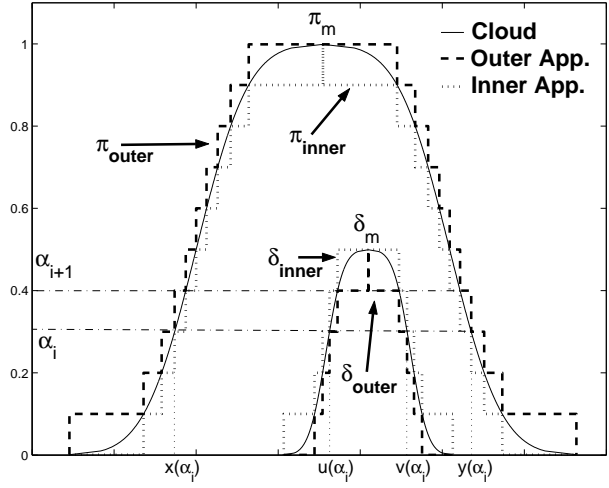}
\caption{Inner and outer approximations of a non-comonotonic clouds}
\label{fig:noncocloudappr}
\end{center}
\end{figure}

Given the above transformations, $\P(\underline{\pi}_n) \subset
\P(\pi) \subset \P(\overline{\pi}_n)$, and  $\lim_{n \to \infty}
\P(\underline{\pi}_n) =  \P(\pi)$ and also $\lim_{n \to \infty}
\P(\overline{\pi}_n) =  \P(\pi)$. Similarly, $\P(1 -
\underline{\delta}_n) \subset \P(1 - \delta) \subset \P(1 -
\overline{\delta}_n)$, $\lim_{n \to \infty}  \P(1 -
\overline{\delta}_n) = \P(1 - \delta)$ and $\lim_{n \to \infty} \P(1
- \underline{\delta}_n) = \P(1 - \delta)$. Since the set of
probabilities induced by the cloud $[\delta, \pi]$ is $\P(\pi) \cap
\P(1 - \delta)$, it is clear that the two credal sets
$\P(\underline{\pi}_n) \cap \P(1 - \underline{\delta}_n) $ and
$\P(\overline{\pi}_n) \cap \P(1 - \overline{\delta}_n)$, are
respectively inner and outer approximations of $\P(\pi) \cap \P(1 -
\delta)$. Moreover:
$$\lim_{n \to \infty}  \P(\underline{\pi}_n) \cap \P(1 - \underline{\delta}_n) =  \P(\pi) \cap \P(1 - \delta)$$   $$\lim_{n \to \infty}  \P(\overline{\pi}_n) \cap \P(1 - \overline{\delta}_n) =  \P(\pi) \cap \P(1 - \delta).$$

The random sets that are inner approximations (by
proposition~\ref{prop:cinnerapp}) of the finite clouds
$[\underline{\delta}_n,\underline{\pi}_n] $ and
$[\underline{\delta}_n,\underline{\pi}_n]$ converge to the
continuous random set defined by the Lebesgue measure on the unit
interval $\alpha \in [0,1]$ and the multimapping $\alpha
\longrightarrow E_{\alpha}$ such that
\begin{displaymath}
E_{\alpha}=\{ r \in \mathbb{R} | (\pi(r) \geq \alpha) \wedge
(\delta(r) < \alpha) \}.
\end{displaymath}
In the limit, it follows that this continuous random set is an inner
approximation of the continuous cloud.
\end{proof}

\end{document}